\documentstyle[12pt,leqno]{article}
\input{mssymb.tex}
\input{epsf.sty}
\font\smit = cmti10 at 11pt

\newcommand{\rme}{{\rm e}}
\newcommand{\exc}{{\rm exc}}
\newcommand{\id}{{\mbox{id}}}
\newcommand{\inv}{{\rm inv}}
\newcommand{\Sing}{{\mbox{Sing}}}
\newcommand{\supp}{{\mbox{supp}}}

\newcommand{\bs}{{\backslash}}
\newcommand{\BC}{{\Bbb C}}
\newcommand{\BF}{{\Bbb F}}
\newcommand{\BK}{{\Bbb K}}
\newcommand{\BN}{{\Bbb N}}
\newcommand{\BP}{{\Bbb P}}
\newcommand{\BQ}{{\Bbb Q}}
\newcommand{\BR}{{\Bbb R}}

\newcommand{\cE}{{\cal E}}
\newcommand{\cF}{{\cal F}}
\newcommand{\cG}{{\cal G}}
\newcommand{\cH}{{\cal H}}
\newcommand{\cI}{{\cal I}}
\newcommand{\cO}{{\cal O}}
\newcommand{\hcO}{{\widehat\cO}}

\newcommand{\th}{\tilde h}
\newcommand{\tsigma}{{\tilde\sigma}}
\newcommand{\tx}{{\tilde x}}
\newcommand{\ty}{{\tilde y}}
\newcommand{\qed}{\hfill$\Box$}

\newcommand{\um}{{\underline m}}
\newcommand{\linee}{\smallskip 
\begin{center}{\underline{\hspace{2truein}}}\end{center}\medskip}

\newcommand{\lra}{{\rightarrow}}

\newcommand{\ex}{\medskip\noindent {\em Example\ }}
\newcommand{\rem}{\medskip\noindent {\em Remark\ }}
\newcommand{\defn}{\medskip\noindent {\em Definition\ }}
\newcommand{\theorem}{\medskip\noindent{\bf Theorem\ }}
\newcommand{\prop}{\medskip\noindent{\bf Proposition\ }}
\newcommand{\alg}{\medskip\noindent{\bf Algorithm.\quad }}
\newcommand{\prf}{\medskip\noindent{\em Proof.\quad }}
             
\begin{document}
\bibliographystyle{plain}

\thispagestyle{empty}
\title{\bf Resolution of Singularities}
\author{
\smallskip
Edward Bierstone and Pierre D. Milman\\
Department of Mathematics\\ 
University of Toronto\\ 
Toronto, Ontario M5S 3G3}
\date{}
\maketitle

%
\setcounter{page}{1}

\section{Introduction}

Resolution of singularities has a long history that goes back
to Newton in the case of plane curves.
For higher-dimensional singular spaces, the problem was formulated
toward the end of the last century, and it was solved in
general, for algebraic varieties defined over fields of characteristic
zero, by Hironaka in his famous paper [H1] of 1964.
([H1] includes the case of real-analytic spaces;
Hironaka's theorem for complex-analytic spaces is proved
in [H2], [AHV1], [AHV2].)
But Hironaka's result is highly non-constructive.
His proof is one of the longest and hardest in mathematics,
and it seems fair to say that only a handful of mathematicians
have fully understood it.
We are not among them!
Resolution of singularities is used in many areas of mathematics,
but even certain aspects of the theorem (for example,
{\em canonicity}; see 1.11 below) have remained unclear.

This article is an exposition of an elementary constructive
proof of canonical resolution of singularities in characteristic zero.
Our proof was sketched in the hypersurface case in [BM4] and is
presented in detail in [BM5].

When we started thinking about the subject almost twenty years
ago, our aim was simply to understand resolution of singularities.
But we soon became convinced that it should be possible to give
simple direct proofs of at least those aspects of the theorem
that are important in analysis.
In 1988, for example, we published a very simple proof that any 
real-analytic variety is the image by a proper analytic mapping 
of a manifold of the same dimension [BM1].
The latter statement is a real version of a local form of resolution
of singularities, called {\em local uniformization}.

It is the idea of [BM1, Section 4] that we have developed (via
[BM2]) to define a new local invariant for desingularization
that is the main subject of this exposition.
Our invariant $\inv_X(a)$ is a finite sequence (of nonnegative
rational numbers and perhaps $\infty$, in the case of a hypersurface),
defined at each point $a$ of our space $X$.
Such sequences can be compared lexicographically.
$\inv_X(\cdot)$ takes only finitely many maximum values
(at least locally), and we get an algorithm for canonical resolution
of singularities by successively blowing up its maximum loci.
Moreover, $\inv_X(\cdot)$ can be described by local computations
that provide equations for the centres of blowing up.

We begin with an example to illustrate the meaning of resolution
of singularities:

\ex {\em 1.1.}\quad
Let $X$ denote the quadratic cone
$x^2-y^2-z^2=0$ in affine 3-space --- the simplest example
of a singular surface.
\bigskip
\begin{center}
{\hskip .3in{\epsfxsize=2.5in
\epsfbox{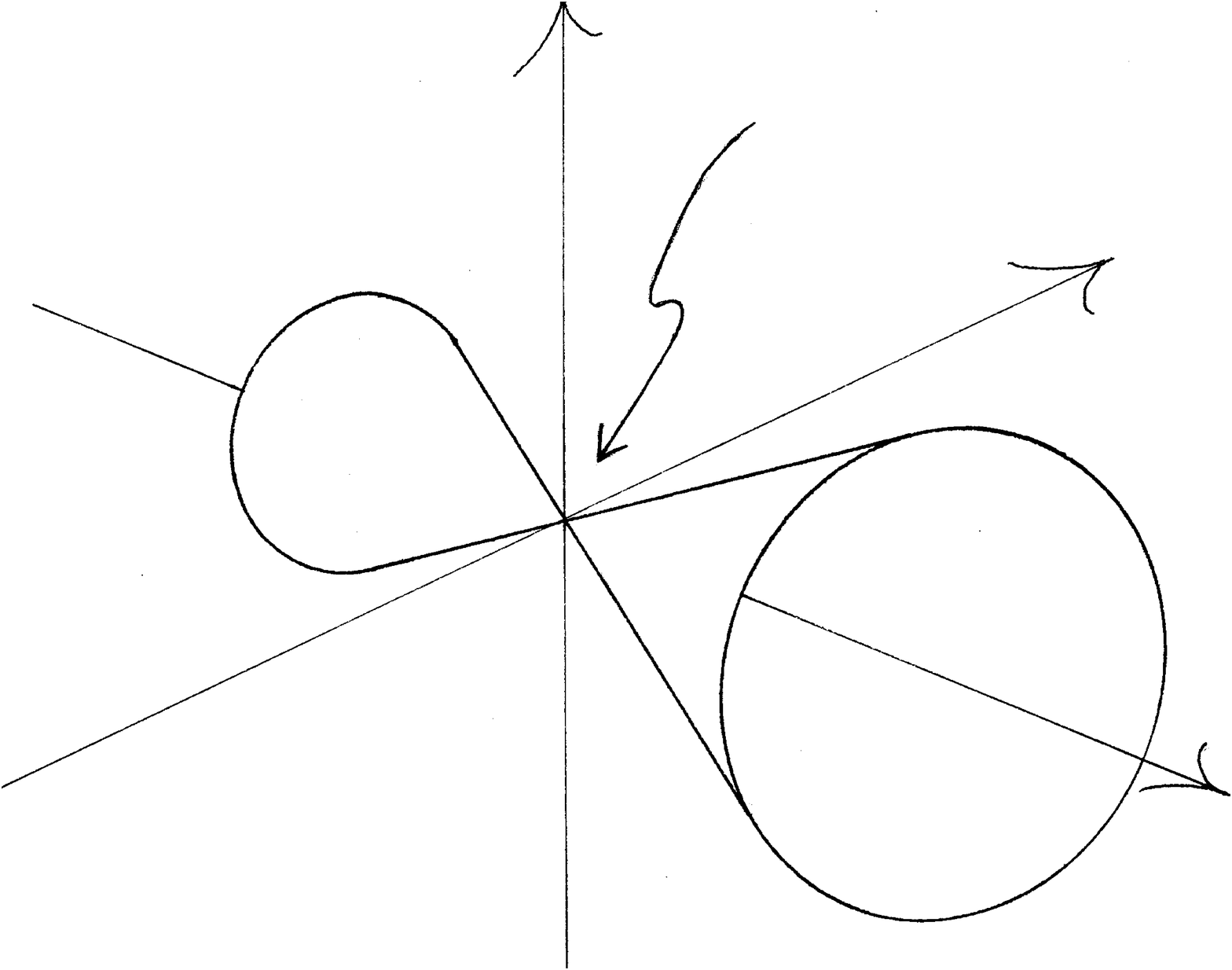}{\vskip -2.18in\hskip .02in{\smit z}}{\vskip .35in \hskip 2.4in{\smit y}}}
{\vskip .95in\hskip 2.8in{\smit x}}{\vskip -1.6in\hskip 1.31in{$\Sing\, X$}}}
\end{center}
\vspace{1.4truein}
\begin{center}
$X:\ x^2-y^2-z^2=0$
\end{center}

$X$ can be desingularized by making a simple quadratic
transformation of the ambient space:
\begin{displaymath}
\sigma:\quad x=u,\ y=uv,\ z=uw .
\end{displaymath}
The inverse image of $X$ by this mapping $\sigma$ is given
by substituting the formulas for $x$, $y$ and $z$ into the
equation of $X$:

\begin{displaymath}
\sigma^{-1}(X):\quad u^2(1-v^2-w^2) = 0.
\end{displaymath}
Thus $\sigma^{-1}(X)$ has two components:
The plane $u=0$ is the set of critical points of the mapping
$\sigma$; it is called the {\em exceptional hypersurface}.
(Here $E':=\{u=0\}$ is the inverse image of the singular
point of $X$.)
The quotient after completely factoring out the ``exceptional
divisor'' $u$ defines what is called the {\em strict
transform\/} $X'$ of $X$ by $\sigma$.
Here $X'$ is the cylinder $v^2+w^2=1$.
\medskip
\begin{center}
{\hskip .1in{\epsfxsize=2.5in 
\epsfbox{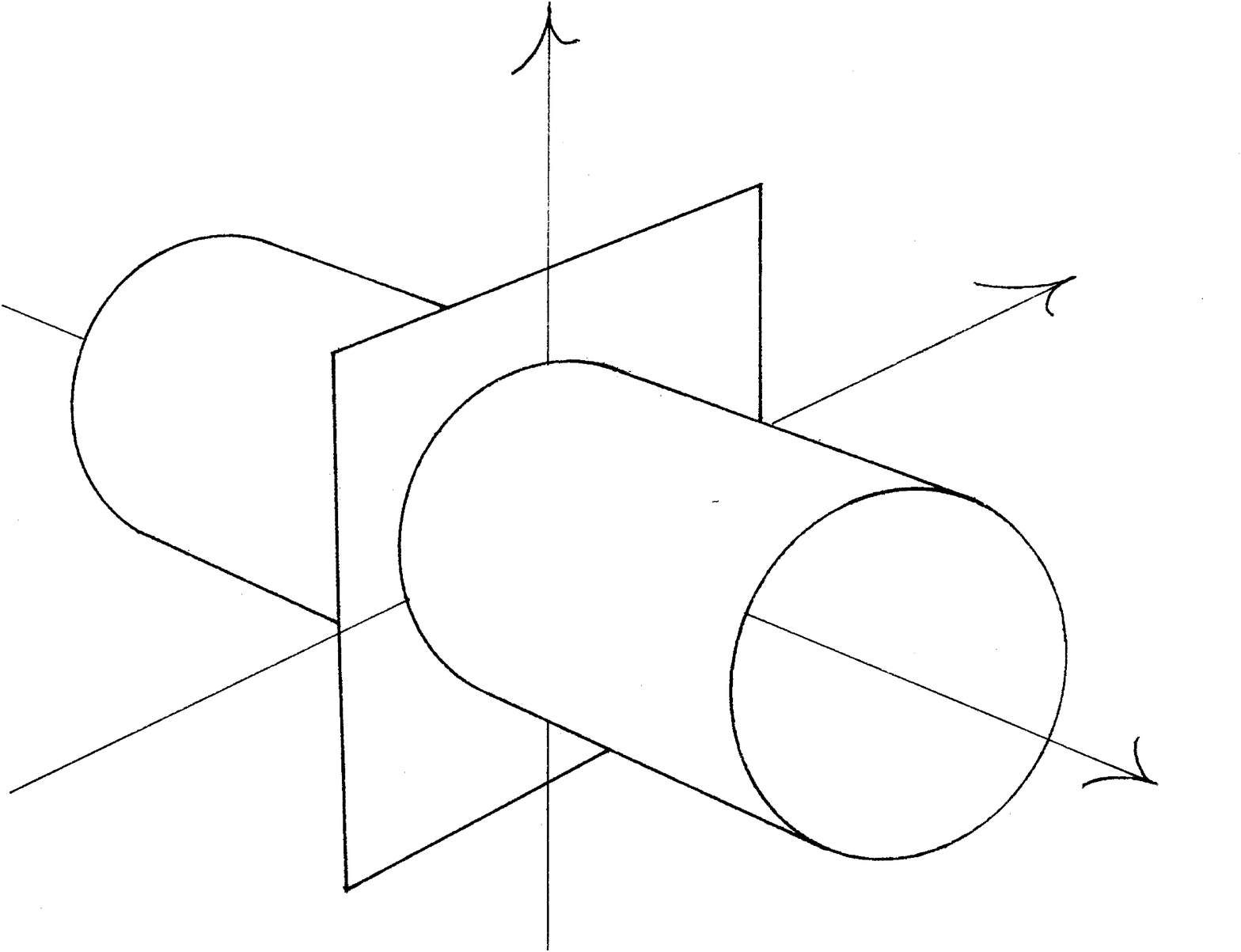}{\vskip -2.18in\hskip -.2in{\smit{w}}}
{\vskip .17in\hskip 1in {$\smit{E':\ u=0}$}}
{\vskip .04in\hskip 2.1in{\smit{v}}}
{\vskip .4in\hskip 3.18in{$\smit{X':\ v^2+w^2=1}$}}
{\vskip .29in\hskip 2.43in{\smit{u}}}}}
\end{center}

\vspace{.4truein}
In this example, $\sigma|X'$ is a 
{\em resolution of singularities\/} of $X$:
$X'$ is smooth and $\sigma|X'$ is a proper mapping
onto $X$ that is an isomorphism outside the singularity.
But the example illustrates a stronger statement, called
{\em embedded resolution of singularities}: $X$ is
desingularized by making a simple transformation
of the ambient space, after which, in addition, the
strict transform $X'$ and the exceptional hypersurface
$E'$ have only {\em normal crossings}; this means
that each point admits a coordinate neighbourhood
with respect to which both $X'$ and $E'$ are coordinate
subspaces.

\medskip

The quadratic transformation $\sigma$ in Example 1.1 is
also called a {\em blowing-up\/} with {\em centre\/} the
origin.
(The centre is the set of critical values of $\sigma$.)
More accurately, the blowing-up of affine 3-space
with centre a point is covered in a natural way
by three affine coordinate charts, and $\sigma$ above is
the formula for the blowing-up restricted to one chart.

Sequences of quadratic transformations, or point blowings-up,
were first used to resolve the singularities of curves
by Max Noether in the 1870's [BN].

The more general statement of ``embedded resolution
of singularities'' seems to have been formulated precisely
first by Hironaka.
But it is implicit already in the earliest rigorous
proofs of local desingularization of surfaces, as a natural
generalization prerequisite to the inductive step of a proof
by induction on dimension (cf. Sections 2,3 below).
For example, in one of the earliest proofs of local
desingularization or uniformization of surfaces, Jung used
embedded desingularization of curves by sequences of
quadratic transformations (applied to the branch locus
of a suitable projection) to prove uniformization
for surfaces [Ju].
Similar ideas were used in the first proofs of global
resolution of singularities of algebraic surfaces,
by Walker [Wa] and Zariski [Z1] in the
late 1930's.
(The latter was the first algebraic proof, by sequences
of normalizations and point blowings-up.)

From the point of view of subsequent work, however,
Zariski's breakthrough came in the early 1940's when
he localized the idea of the centre of blowing-up, thus
making possible an extension of the notion of quadratic
transformation to blowings-up with centres that are
not necessarily $0$-dimensional [Z2].
This led Zariski to a version of embedded resolution
of singularities of surfaces, and to a weaker (non-embedded)
theorem for 3-dimensional algebraic varieties [Z3].
It was the path that led to Hironaka's great theorem
and to most subsequent work in the area, including our own.
(See References below.)

From a general viewpoint, some important features of our work
in comparison with previous treatments are: (1) It is
canonical. (See 1.11.)
(2) We isolate simple local properties of an invariant
(Section 4, Theorem B) from which global desingularization is automatic.
(3) Our proof in the case of a hypersurface (a space defined
locally by a single equation) does not involve passing
to higher codimension (as in the inductive procedure of
[H1]).

Very significant results on resolution of singularities
over fields of nonzero characteristic have recently been
obtained by de Jong [dJ] and have been announced
by Spivakovsky.

\medskip

\noindent {\em 1.2. Blowing up.} We first describe the blowing-up
of an open subset $W$ of $r$-dimensional affine space
with centre a point $a$.
(Say $a=0\in W$.)
The {\em blowing-up\/} $\sigma$ with {\em centre\/} $0$
is the projection onto $W$ of a space $W'$ that is
obtained by replacing the origin by the $(r-1)$-dimensional
projective space $\BP^{r-1}$ of all lines through $0$:
\[
W' = \{\, (x,\lambda)\in W\times \BP^{r-1}:\ x\in\lambda\,\}
\]
and $\sigma$: $W'\rightarrow W$ is defined by $\sigma(x,\lambda)=x$.
(Outside the origin, a point $x$ belongs to a unique line
$\lambda$, but $\sigma^{-1}(0) = \BP^{r-1}$.
Clearly, $\sigma$ is a proper mapping.)
$W'$ has a natural algebraic structure:
If we write $x$ in terms of the affine coordinates
$x=(x_1,\ldots,x_r)$, and $\lambda$ in the corresponding
homogeneous coordinates $\lambda=[\lambda_1,\ldots,\lambda_r]$,
then the relation $x\in\lambda$ translates into the system
of equations $x_i\lambda_j = x_j\lambda_i$, for all $i,j$.

These equations can be used to see that $W'$ has the
structure of an algebraic manifold:
For each $i=1,\ldots,r$, let $W_i'$ denote the open
subset of $W'$ where $\lambda_i\ne 0$.
In $W_i'$, $x_j=x_i\lambda_j/\lambda_i$, for each $j\ne i$,
so we see that $W_i'$ is smooth:
it is the graph of a mapping in terms of coordinates
$(y_1,\ldots,y_r)$ for $W_i'$ defined by
$y_i=x_i$ and $y_j=\lambda_j/\lambda_i$ if $j\neq i$.
In these coordinates, $\sigma$ is a quadratic transformation
given by the formulas
\[
x_i=y_i,\quad x_j=y_iy_j\ \  \mbox{for all $j\ne i$,}
\]
as in Example 1.1.

Once blowing up with centre a point has been described
as above, it is a simple matter to extend the idea to
blowing up a manifold, or smooth space, $M$ with centre
an arbitrary smooth closed subspace $C$ of $M$:
Each point of $C$ has a product coordinate neighbourhood
$V\times W$ in which $C = V\times\{0\}$; over this
neighbourhood, the blowing-up with centre $C$ identifies
with $\id_V\times\sigma$:
$V\times W'\rightarrow V\times W$, where $\id_V$ is the
identity mapping of $V$ and $\sigma$: $W'\rightarrow W$
is the blowing-up of $W$ with centre $\{0\}$.
The blowing-up $M'\rightarrow M$ with centre $C$ is an
isomorphism over $M\bs C$.
The preceding conditions determine $M'\rightarrow M$ uniquely,
up to an isomorphism of $M'$ commuting with the projections
to $M$.

\ex {\em 1.3.}\quad 
\begin{center}
{\hskip .1in{\epsfxsize=2.5in
\epsfbox{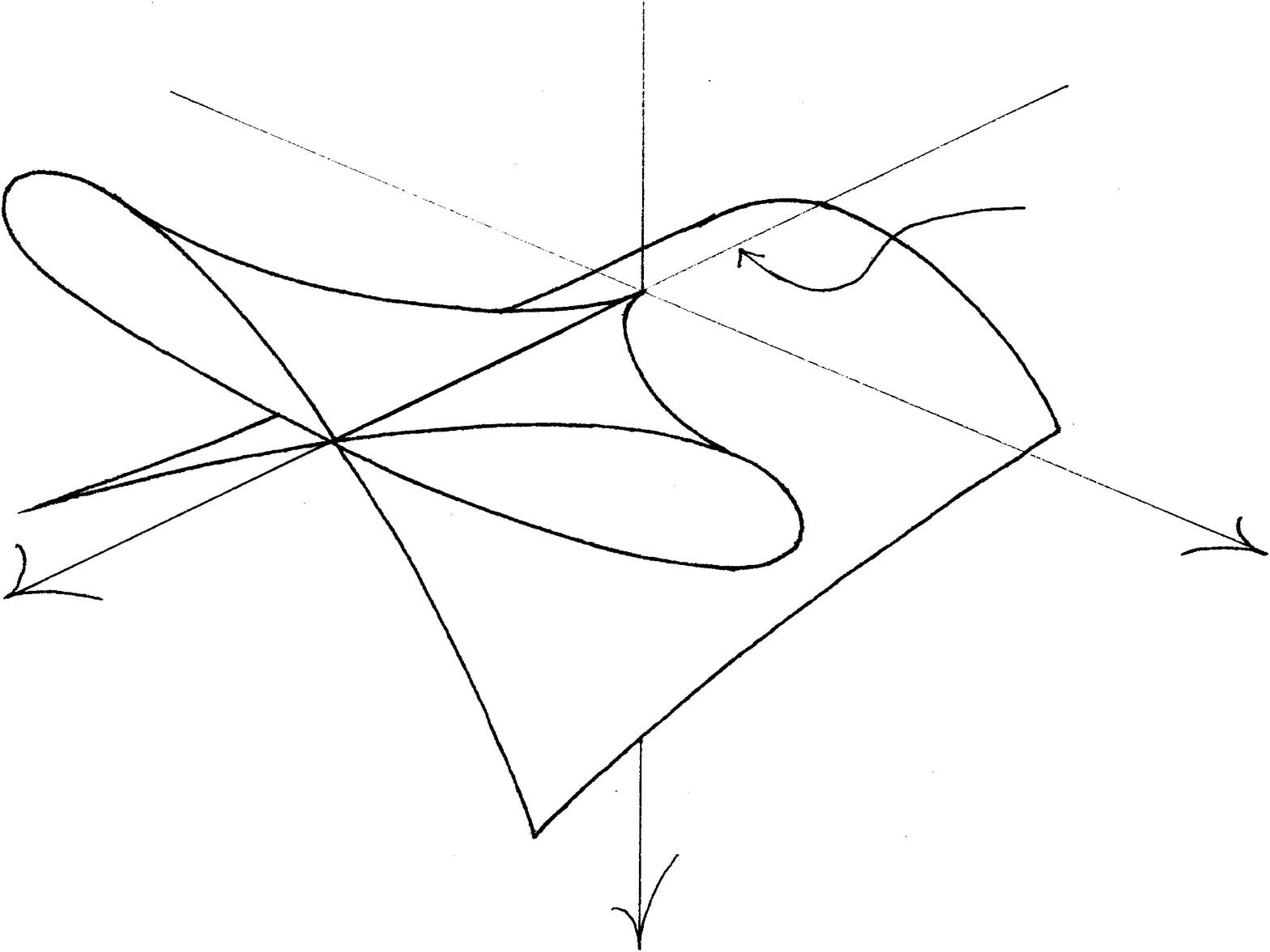}{\vskip -1.59in\hskip 2.19in{$\Sing\, X$}}}
{\vskip .45in\hskip 2.6in{\smit x}}
{\vskip 0.6in \hskip .11in{\smit z}}
{\vskip -.91in \hskip -2.4in{\smit y}}}
\end{center}

\vspace{.5truein}

\begin{center}
{$X:\ z^3-x^2yz-x^4=0$}
\end{center}
\vspace{.5truein}

This surface is particularly interesting in the real
case because, as a {\em subset\/} of $\BR^3$, it is
singular only along the nonnegative $y$-axis.
But resolution of singularities is an algebraic process:
it applies to spaces that include a functional structure
(given here by the equation for $X$).
As a {\em subspace\/} of $\BR^3$, $X$ is singular along
the entire $y$-axis.

\medskip

In general for a hypersurface $X$ ---  say that $X$ is
defined locally by an equation $f(x)=0$ --- to say that
a point $a$ is {\em singular\/} means there are no
linear terms in the Taylor expansion of $f$ at $a$;
in other words, the order $\mu_a(f) > 1$.
(The {\em order\/} or {\em multiplicity\/} $\mu_a(f)$
of $f$ at $a$ is the degree of the lowest-order
homogeneous part of the Taylor expansion of $f$ at $a$.
We will also call $\mu_a(f)$ the {\em order\/} $\nu_{X,a}$
of the hypersurface $X$ at $a$.)

\medskip

The general philosophy of our approach to desingularization
(going back to Zariski [Z3]) is the blow up with smooth
centre as large as possible inside the locus of the most
singular points.
In our example here, $X$ has order $3$ at each point of the
$y$-axis.
In general, order is not a delicate enough invariant
to determine a centre of blowing-up for resolution
of singularities, even in the hypersurface case.
(We will refine order in our definition of $\inv_X$.)
But here let us take the blowing-up $\sigma$ with centre
the $y$-axis:
\[
\sigma:\quad x=u,\ y=v,\ z=uw.
\]
(Again, this is the formula for blowing up in one of
two coordinate charts required to cover our space.
But the strict transform of $X$ in fact lies completely
within this chart.)
The inverse image of $X$ is
\[
\sigma^{-1}(X):\quad u^3(w^3-vw-u)=0;
\]
$\{u=0\}$ is the exceptional hypersurface $E'$ (the
inverse image of the centre of blowing up) and the
strict transform $X'$ is smooth.
(It is the graph of a function $u=w^3-vw$.)

\begin{center}
{\hskip .1in{\epsfxsize=2.5in
\epsfbox{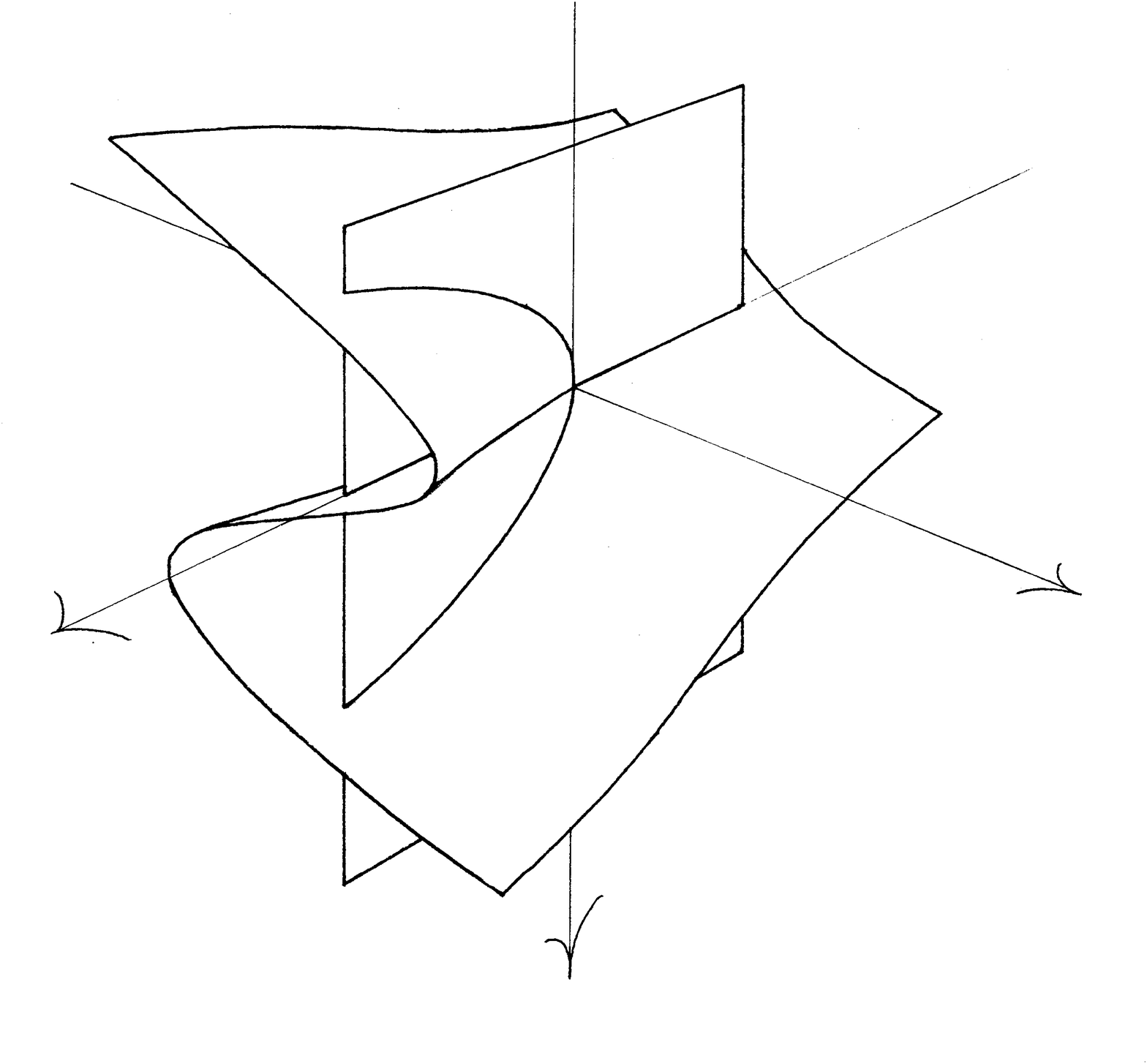}}{\vskip -1.08in\hskip -2.08in{\smit v}}
{\vskip .55in\hskip .18in{\smit w}}
{\vskip -1.04in\hskip 2.39in{\smit u}}
{\vskip -.6in\hskip 3.05in{$X':\ u=w^3-vw$}}
{\vskip -.9in\hskip 1.7in{$E':\ u=0$}}}
\end{center}

\vspace{2truein}

$X'$ is a desingularization of $X$, but we have not yet
achieved an embedded resolution of singularities because
$X'$ and $E'$ do not have normal crossings at the origin.
Further blowings-up are needed for embedded resolution
of singularities.
\medskip

\noindent {\em 1.4. Embedded resolution of singularities.}
Let $X$ denote a (singular) space.
We assume, for simplicity, that $X$ is a closed subspace
of a smooth ambient space $M$.
(This is always true locally.)
The goal of embedded desingularization, in its simplest
version, is to find a proper morphism $\sigma$ from a
smooth space $M'$ onto $M$, in our category, with
the following properties:

(1) $\sigma$ is an isomorphism outside the singular
locus $\Sing\, X$ of $X$.

(2) The strict transform $X'$ of $X$ by $\sigma$ is
smooth.
(See 1.6 below.)
$X'$ can be described geometrically (at least if our
field $\BK$ is algebraically closed;
cf. [BM5, Rmk. 3.15]) as the smallest closed subspace
of $M'$ that includes $\sigma^{-1}(X\bs\Sing\, X)$.

(3) $X'$ and $E'=\sigma^{-1}(\Sing\, X)$ simultaneously
have only normal crossings.
This means that, locally, we can choose coordinates
with respect to which $X'$ is a coordinate subspace
and $E'$ is a collection of coordinate hyperplanes.

We can achieve this goal with $\sigma$ the composite
of a sequence of blowings-up; a finite sequence when our
spaces have a compact topology (for example, in an algebraic
category), or a locally-finite sequence for non-compact
analytic spaces.
(A sequence of blowings-up over $M$ is {\em locally finite\/}
if all but finitely many of the blowings-up are trivial
over any compact subset of $M$.
The composite of a locally-finite sequence of blowings-up
is a well-defined morphism $\sigma$.)

\medskip

\noindent {\em 1.5. The category of spaces.}
Our desingularization theorem applies to the usual spaces
of algebraic and analytic geometry over fields $\BK$ of
characteristic zero --- algebraic varieties, schemes
of finite type, analytic spaces (over $\BR$, $\BC$ or any
locally compact $\BK$) --- but in addition to certain
categories of spaces intermediate between analytic
and $C^\infty$ (See [BM5].)
In any case, we are dealing with a category of local-ringed
spaces $X=(|X|,\cO_X)$ over $\BK$, where $\cO_X$ is a
coherent sheaf of rings.
We are intentionally not specific about the category
in this exposition because we want to emphasize
the principles involved, and the main requirement
for our desingularization algorithm is simply that
a smooth space $M=(|M|,\cO_M)$ in our category admit
a covering by {\em (regular) coordinate charts\/} in which
we have analogues of the usual operations of calculus
of analytic functions; namely:

The coordinates $(x_1,\ldots,x_n)$ of a chart $U$ are
{\em regular functions\/} on $U$ (i.e., each $x_i\in
\cO_M(U)$) and all partial derivatives $\partial^{|\alpha|}/
\partial x^\alpha = \partial^{\alpha_1+\cdots+\alpha_n}/
\partial x_1^{\alpha_1}\cdots \partial x_n^{\alpha_n}$
make sense as transformations $\cO_M(U)\rightarrow \cO_M(U)$.
Moreover, for each $a\in U$, there is an {\em injective\/} ``Taylor
series homomorphism'' $T_a$: $\cO_{M,a}\rightarrow \BF_a[[X]]=
\BF_a[[X_1,\ldots,X_n]]$, where $\BF_a$ denotes the
residue field $\cO_{M,a}/\um_{M,a}$, such that $T_a$
induces an isomorphism $\hcO_{M,a}\stackrel{\cong}{\lra}
\BF_a [[X]]$ and $T_a$
commutes with differentiation: $T_a\circ (\partial^{|\alpha|}/
\partial x^\alpha)=(\partial^{|\alpha|}/\partial X^\alpha)\circ
T_a$, for all $\alpha\in\BN^n$.
($\um_{M,a}$ denotes the maximal ideal and $\hcO_{M,a}$
the completion of $\cO_{M,a}$.
$\BN$ denotes the
nonnegative integers.)

\medskip

In the case of real- or complex-analytic spaces, of
course, $\BK=\BR$ or $\BC$, $\BF_a=\BK$ at each point,
and ``coordinate chart'' means the classical notion.
Regular coordinate charts for schemes of finite type
are introduced in [BM5, Section 3].

Suppose that $M=(|M|,\cO_M)$ is a manifold (smooth space)
and that $X=(|X|,\cO_X)$ is a closed subspace of $M$.
This means there is a coherent sheaf of ideals $\cI_X$
in $\cO_M$ such that $|X|=\supp\,\cO_M/\cI_X$ and
$\cO_X$ is the restriction to $|X|$ of $\cO_M/\cI_X$.
We say that $X$ is a {\em hypersurface\/} in $M$ if
$\cI_{X,a}$ is a principal ideal, for each $a\in |X|$.
Equivalently, for every $a\in |X|$, there is an open
neighbourhood $U$ of $a$ in $|M|$ and a regular function
$f\in\cO_M(U)$ such that $|X|U|=\{\, x\in U:\ f(x)=0\,\}$
and $\cI_X|U$ is the principal ideal $(f)$ generated
by $f$; we write $X|U = V(f)$.

\medskip

\noindent {\em 1.6. Strict transform.}
Let $X$ denote a closed subspace of a manifold $M$, and let $\sigma$:
$M'\rightarrow M$ be a blowing-up with smooth centre $C$.
If $X$ is a hypersurface, then the strict transform $X'$
of $X$ by $\sigma$ is a closed subspace of $M'$ that can be
defined as follows:
Say that $X=V(f)$ in a neighbourhood of $a\in|X|$.
Then, in some neighbourhood of $a'\in \sigma^{-1} (a)$,
$X'=V(f')$, where $f'=y_\exc^{-d} f\circ \sigma$, $y_\exc$
denotes a local generator of $\cI_{\sigma^{-1}(C)}\subset\cO_{M'}$,
and $d=\mu_{C,a}(f)$ denotes the {\em order\/} of $f$ {\em along\/}
$C$ at $a$:
$d=\max\{\,k:\ (f)\subset \cI_{C,a}^k\,\}$; $d$ is the largest
power to which $y_\exc$ factors from $f\circ\sigma$ at $a'$.

The strict transform $X'$ of a general closed subspace $X$ of
$M$ can be defined locally, at each $a'\in \sigma^{-1}(a)$,
as the intersection of all hypersurfaces $V(f')$, for all
$f\in\cI_{X,a}$.
We likewise define the strict transform by a sequence of blowings-up
with smooth centres.

Each of the categories listed in 1.5 above is closed under
blowing up and strict transform [BM5, Prop. 3.13 ff.];
the latter condition is needed to apply the desingularization
algorithm in a given category.

\medskip

\noindent {\em 1.7. The invariant.}
Let $X$ denote a closed subspace of a manifold $M$.
To describe $\inv_X$, we consider a sequence of transformations
$$ 
\begin{array}{cccccccccc}
\lra & M_{j+1} & \stackrel{\sigma_{j+1}}{\lra}
& M_j & \lra & \cdots & \lra & M_1 & \stackrel{\sigma_1}{\lra}
& M_0=M \\
& X_{j+1} & & X_j & & & & X_1 & & X_0 = X \\
& E_{j+1} & & E_j & & & & E_1 & & E_0 = \emptyset 
\end{array}
\leqno(1.8)
$$ 
where, for each $j$, $\sigma_{j+1}$: $M_{j+1}\rightarrow M_j$
denotes a blowing-up with smooth centre $C_j\subset M_j$,
$X_{j+1}$ is the strict transform of $X_j$ by $\sigma_{j+1}$,
and $E_{j+1}$ is the set of exceptional hypersurfaces in $M_{j+1}$;
i.e., $E_{j+1}=E_j'\cup \{\sigma_{j+1}^{-1} (C_j)\}$, where $E_j'$
denotes the set of strict transforms by $\sigma_{j+1}$ of all
hypersurfaces in $E_j$.

Our invariant $\inv_X(a)$, $a\in M_j$, $j=0,1,2,\ldots$, will
be defined inductively over the sequence of blowings-up;
for each $j$, $\inv_X(a)$, $a\in M_j$, can be defined provided
that the centres $C_i$, $i<j$, are {\em admissible\/} (or
$\inv_X$-{\em admissible}) in the sense that:

(1) $C_i$ and $E_i$ simultaneously have only normal crossings.

(2) $\inv_X(\cdot)$ is locally constant on $C_i$.

The condition (1) guarantees that $E_{i+1}$ is a collection of smooth
hypersurfaces having only normal crossings.
We can think of the desingularization algorithm in the following way:
$X\subset M$ determines $\inv_X(a)$, $a\in M$, and thus the
first admissible centre of blowing up $C=C_0$; then $\inv_X(a)$
can be defined on $M_1$ and determines $C_1$, etc.

The notation $\inv_X(a)$, where $a\in M_j$, indicates a
dependence not only on $X_j$, but also on the original space
$X$.
In fact $\inv_X(a)$, $a\in M_j$, is invariant under local
isomorphisms of $X_j$ that preserve $E(a)=\{\, H\in E_j:\ H\ni a\,\}$
and certain subcollections $E^r(a)$ (which will be taken to
encode the history of the resolution process).
To understand why some dependence on the history should be needed,
let us consider how, in principle, it might be possible to determine
a {\em global\/} centre of blowing up using a {\em local\/} invariant:

\ex {\em 1.9.}\quad 
It is easy to find an example of a surface $X$ whose singular locus,
in a neighbourhood of a point $a$, consists of two smooth curves
\begin{minipage}[t]{2in}
with a normal crossing at $a$, and where $X$ has the property
that, if we blow up with centre $\{ a\}$, then there are
points $a'$ in the fibre $\sigma^{-1}(a)$ where the strict
transform $X'$ has the {\em same\/} local equation (in suitable\break 
\vspace{-.15in}
\end{minipage}\hfill
\raisebox{.4cm}{\begin{minipage}[t]{3in}
\hskip 4in{\epsfxsize=2.5in \epsfbox{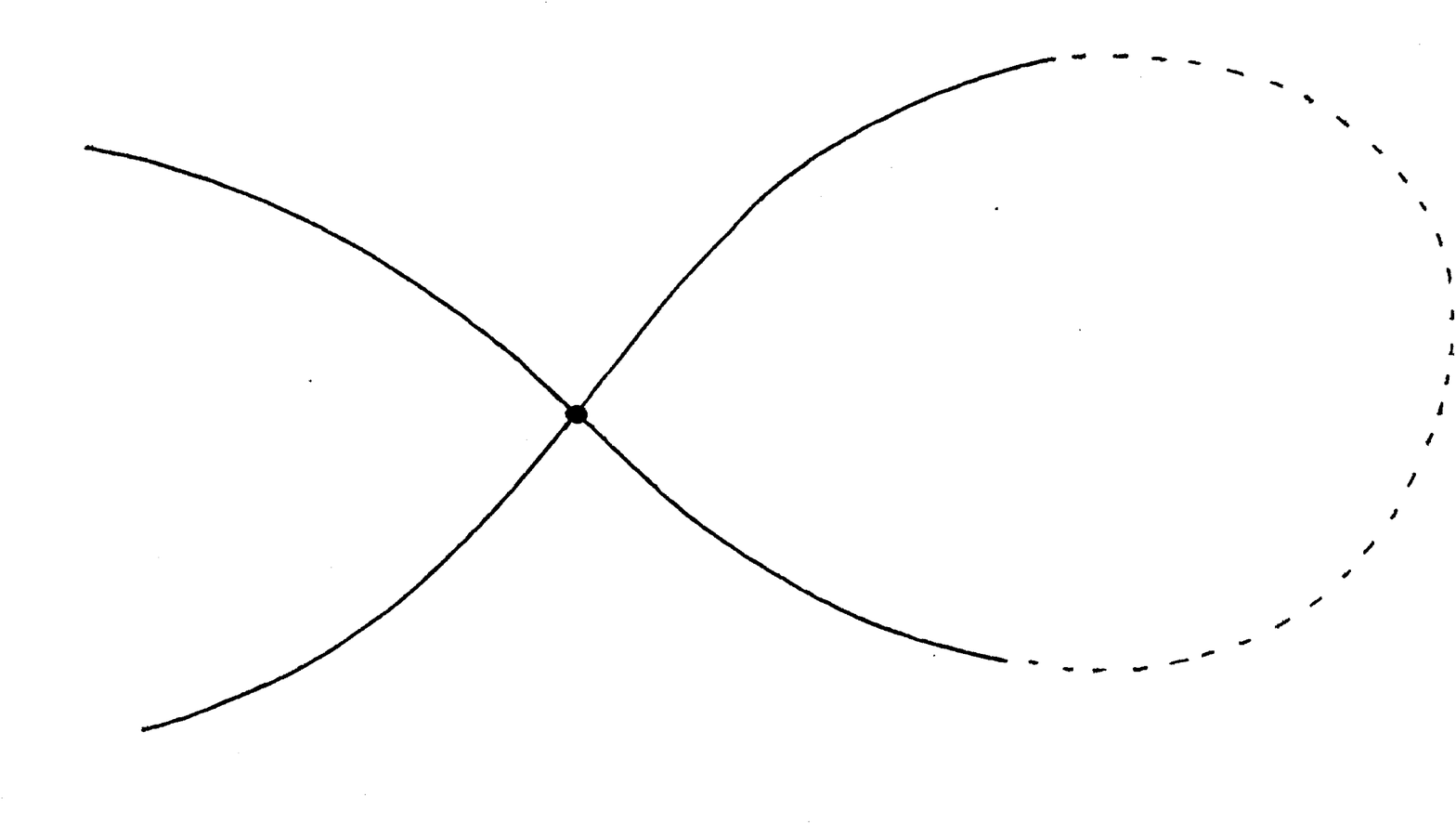}{\vskip -1in\hskip .87in{$a$}}{\vskip .5in\hskip .3in{$\Sing\, X$}}}
\end{minipage}}
coordinates) as that of $X$ at $a$, or an even more complicated
equation (as in Example 3.1 below). 
This suggests that to simplify the singularities in a neighbourhood
of $a$ by blowing up with smooth centre in $\Sing\, X$, we should
choose as centre one of the two smooth curves.
But our surface may have the property that neither curve extends
to a global smooth centre, as illustrated.
So there is no choice but to blow up with centre $\{ a\}$, although
it seems to accomplish nothing: The figure shows the \break
{\begin{minipage}[t]{2in}
singular locus of $X'$; there are two points
$a'\in \sigma^{-1} (a)$ where the singularity is the
same as or worse than before.
But what has changed at each of these points is the status of one
of the curves, which is now {\em exceptional}.
The moral is that,\break 
\vspace{-.13in}
\end{minipage}}\hfill
\raisebox{.1cm}{\begin{minipage}[t]{3in}
\hskip 4in{\epsfxsize=2.5in\epsfbox{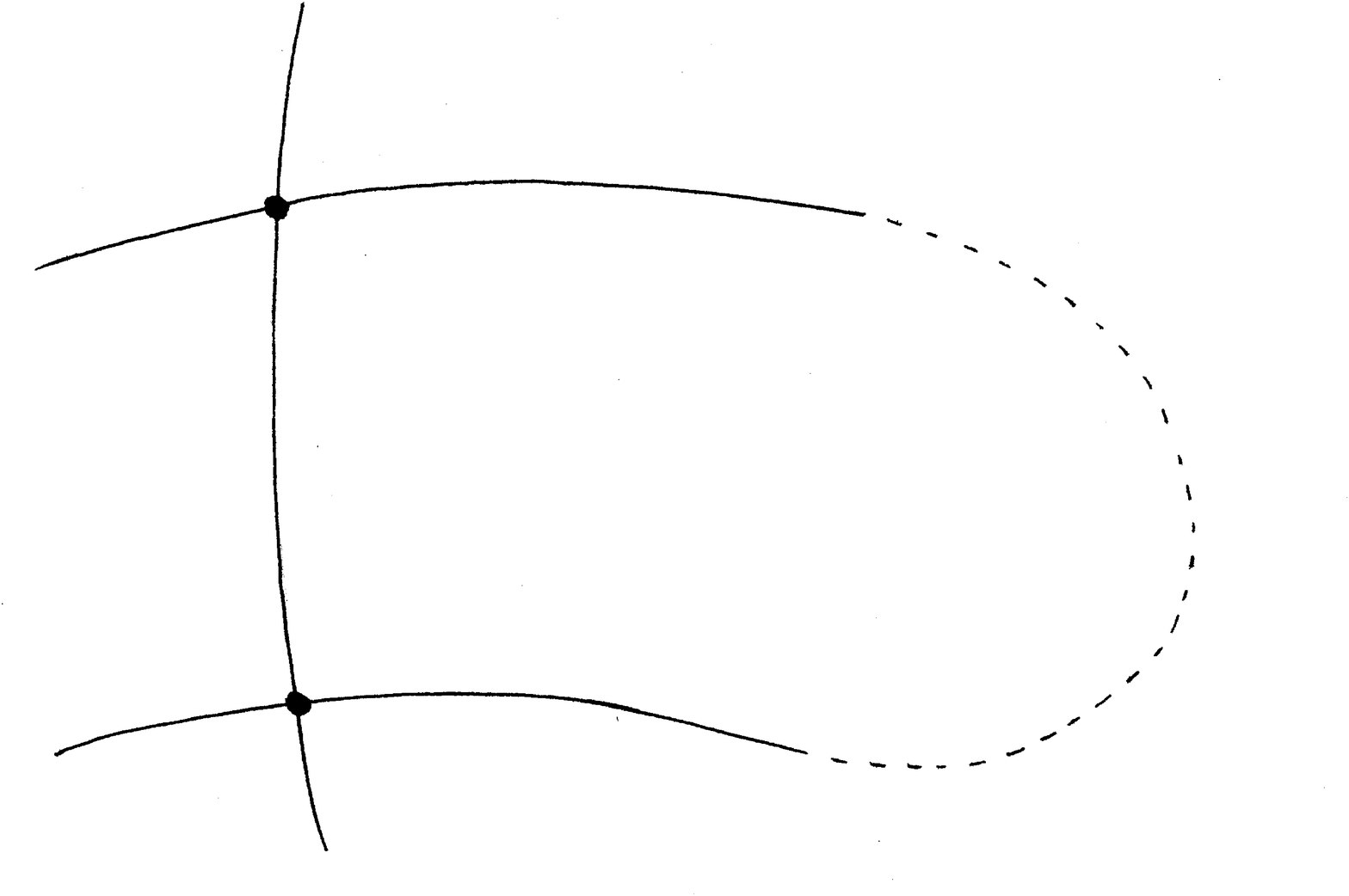}{\vskip -1.7in\hskip .62in{$X'\cap E'$}}{\vskip .86in\hskip .61in{$a'$}}}
\end{minipage}}
\noindent although the singularity of $X$ at $a$ has not
been simplified in the strict transform, an invariant which takes
into account the history of the resolution process as recorded
by the accumulating exceptional hypersurfaces might nevertheless
measure some improvement.
\medskip

Consider a sequence of blowings-up as before.
For simplicity, {\em we will assume that\/} $X\subset M$ 
{\em is a hypersurface}.
Then $\inv_X(a)$, $a\in M_j$, is a finite sequence beginning
with the order $\nu_1(a)=\nu_{X_j,a}$ of $X_j$ at $a$:
\[
\inv_X(a)\ =\ \big( \nu_1(a), s_1(a);\, \nu_2(a),s_2(a);\,
\ldots,\, s_t(a);\, \nu_{t+1}(a)\big) .
\]
(In the general case, $\nu_1(a)$ is replaced by a more delicate
invariant of $X_j$ at $a$ --- the Hilbert-Samuel function
$H_{X_j,a}$ (see [BM5]) --- but the remaining entries
of $\inv_X(a)$ are still rational numbers (or $\infty$) as we will
describe, and the theorems below are unchanged.)
The $s_r(a)$ are nonnegative integers counting exceptional
hypersurfaces that accumulate in certain blocks $E^r(a)$ depending
on the history of the resolution process.
And the $\nu_r(a)$, $r\ge 2$, represent certain ``higher-order
multiplicities'' of the equation of $X_j$ at $a$;
$\nu_2(a),\ldots,\nu_t(a)$ are quotients of positive integers
whose denominators are bounded in terms of the previous entries
of $\inv_X(a)$.
(More precisely, $e_{r-1}! \nu_r (a)\in \BN$, $r=1,\ldots,t$,
where $e_0=1$ and $e_r = \max\{ e_{r-1}!, e_{r-1}! \nu_r
(a)\}$.)
The pairs $\big( \nu_r(a), s_r(a)\big)$ can be defined successively
using data that depends on $n-r+1$ variables (where $n$ is the ambient
dimension), so that $t\le n$ by exhaustion of variables; the final
entry $\nu_{t+1}(a)$ is either $0$ (the order of a nonvanishing
function) or $\infty$ (the order of the function identically zero).

\ex {\em 1.10.}\quad 
Let $X\subset\BK^n$ be the hypersurface
$x_1^{d_1} + x_2^{d_2}+\cdots + x_t^{d_t} =0$, where $1<d_1\le \cdots
\le d_t$, $t\le n$.
Then
\[
\inv_X(0) = \left( d_1, 0;\, \frac {d_2}{d_1}, 0;\, \ldots;\,
\frac{d_t}{d_{t-1}},0;\, \infty\right) .
\]
This is $\inv_X(0)$ in ``year zero'' (before the first blowing up),
so there are no exceptional hypersurfaces.

\theorem {\bf A.}\quad {\rm
(Embedded desingularization.)}
{\em There is a finite sequence of blowings-up (1.8) with smooth
$\inv_X$-admissible centres $C_j$ (or a locally finite sequence,
in the case of noncompact analytic spaces) such that:

{\rm (1)} For each $j$, either $C_j\subset\Sing\, X_j$ or $X_j$
is smooth and $C_j\subset X_j\cap E_j$.

{\rm (2)} Let $X'$ and $E'$ denote the final strict transform
of $X$ and exceptional set, respectively.
Then $X'$ is smooth and $X'$, $E'$ simultaneously have only normal
crossings.}
\medskip

If $\sigma$ denotes the composite of the sequence of blowings-up
$\sigma_j$, then $E'$ is the critical locus of $\sigma$ and 
$E'=\sigma^{-1} (\Sing\, X)$.
In each of our categories of spaces, $\Sing\, X$ is closed in
the {\em Zariski topology\/} of $|X|$ (the topology whose closed
sets are of the form $|Y|$, for any closed subspace $Y$ of $X$;
see [BM5, Prop. 10.1]).
Theorem A resolves the singularities of $X$ in a meaningful geometric
sense provided that $|X|\bs \Sing\, X$ is (Zariski-)dense in $|X|$.
(For example, if $X$ is a {\em reduced\/} complex-analytic space
or a scheme of finite type.)
More precise desingularization theorems (for example, for spaces
that are not necessarily reduced) are given in [BM5, Ch. IV].

This paper contains an essentially complete proof of Theorem A in
the hypersurface case, presented though in a more informal way than
in [BM5].
We give a constructive definition of $\inv_X$ in Section 3, in
parallel with a detailed example.
In Section 4, we show that $\inv_X$ is indeed an invariant,
and we summarize its key properties in Theorem B.
(The terms $s_r(a)$ of $\inv_X(a)$ can, in fact, be introduced
immediately in an invariant way; see 1.12 below.)
It follows from Theorem B(3) that the maximum locus of $\inv_X$
has only normal crossings and, moreover, each of its local
components extends to a global smooth subspace.
(See Remark 3.6.)

The point is that each component is the intersection of the 
maximum locus of $\inv_X$ with those exceptional hypersurfaces
containing the component;
the exceptional divisors serve as global coordinates.)
We can obtain Theorem A by successively blowing up with centre
given by any component of the maximum locus.

\medskip

\noindent {\em 1.11. Universal and canonical desingularization.}
The exceptional hypersurfaces (the elements of $E_j$) can be ordered
in a natural way (by their ``years of birth'' in the history
of the resolution process).
We can use this ordering to extend $\inv_X(a)$ by an additional
term $J(a)$ that will have the effect of picking out one component
of the maximum locus of $\inv_X(\cdot)$ in a canonical way;
see Remark 3.6.
We write $\inv_X^\rme (\cdot)$ for the extended invariant
$\big( \inv_X(\cdot); J(\cdot)\big)$.
Then our embedded desingularization theorem A can be obtained
by the following:

\alg
{\em Choose as each successive centre of blowing
up $C_j$ the maximum locus of $\inv_X^\rme$ on $X_j$.}
\medskip

The algorithm stops when our space is ``resolved'' as in
the conclusion of Theorem A.
In the general (not necessarily hypersurface) case,
we choose more precisely as each successive centre
$C_j$ the maximum locus of $\inv_X^\rme$ on the non-resolved
locus $Z_j$ of $X_j$; in general, $\{ x:\ \inv_X(x)=\inv_X(a)\}
\subset Z_j$ (as germs at $a$), so that again each $C_j$
is smooth, by Theorem B(3), and the algorithm stops when
$Z_j=\emptyset$.

The algorithm applies to a category of spaces satisfying
a compactness assumption (for example, schemes of finite type,
restrictions of analytic spaces to relatively compact open
subsets), so that $\inv_X(\cdot)$ has global maxima.
Since the centres of blowing up are completely determined by an
invariant, our desingularization theorem is automatically
{\em universal\/} in the following sense:
To every $X$, we associate a morphism of resolution of singularities
$\sigma_X$: $X'\rightarrow X$ such that any local isomorphism 
$X|U \rightarrow Y|V$ (over open subsets $U$ of $|X|$ 
and $V$ of $|Y|$) lifts to an isomorphism $X'|\sigma_X^{-1} 
(U)\rightarrow Y'|\sigma_Y^{-1} (V)$
(in fact, lifts to isomorphisms throughout the entire towers
of blowings-up).

For analytic spaces that are not necessarily compact, we can use
an exhaustion by relatively compact open sets to deduce
{\em canonical\/} resolution of singularities:
Given $X$, there is a morphism of desingularization
$\sigma_X$: $X'\rightarrow X$ such that any local isomorphism
$X|U \rightarrow X|V$ (over open subsets of $|X|$) lifts to
an isomorphism $X'|\sigma_X^{-1} (U) \rightarrow X'|
\sigma_X^{-1} (V)$.  (See [BM5, Section 13].)

\medskip

\noindent {\em 1.12. The terms $s_r(a)$.}
The entries $s_1(a),\, \nu_2(a),\, s_2(a),\,\ldots$ of $\inv_X(a) =
\big( \nu_1(a), s_1(a);\, \ldots, s_t(a);\, \nu_{t+1}(a)\big)$
will themselves be defined recursively.
Let us write $\inv_r$ for $\inv_X$ truncated after $s_r$
(with the convention that $\inv_r(a) = \inv_X(a)$ if $r>t$).
We also write $\inv_{r+\frac{1}{2}}=(\inv_r;\nu_{r+1})$ (with
the same convention), so that $\inv_{\frac{1}{2}}(a)$ means
$\nu_1(a)=\nu_{X_j,a}$ (in the hypersurface case, or
$H_{X_j,a}$ in general).
For each $r$, the entries $s_r$, $\nu_{r+1}$ of $\inv_X$ can be
defined over a sequence of blowings-up (1.8) whose centres
$C_i$ are $(r-\frac{1}{2})$-{\em admissible\/} 
(or $\inv_{r-\frac{1}{2}}$-{\em admissible}) in the sense that:

(1) $C_i$ and $E_i$ simultaneously have only normal crossings.

(2) $\inv_{r-\frac{1}{2}}(\cdot)$ is locally constant on $C_i$.

The terms $s_r(a)$ can be introduced immediately, as follows:
Write $\pi_{ij}=\sigma_{i+1}\circ\cdots\cdot \sigma_j$, $i=0,\ldots,
j-1$, and $\pi_{jj}=$ identity.
If $a\in M_j$, set $a_i=\pi_{ij}(a)$, $i=0,\ldots,j$.
First consider a sequence of blowings-up (1.8) with
$\frac{1}{2}$-admissible centres.
($\inv_{\frac{1}{2}}=\nu_1$ can only decrease over such a sequence;
see, for example, Section 2 following.)
Suppose $a\in M_j$.
Let $i$ denote the ``earliest year'' $k$ such that 
$\nu_1(a)=\nu_1(a_k)$, and set $E^1(a)=\{\, H\in E(a):$ $H$ is 
the strict transform of some hypersurface in $E(a_i)\,\}$.
We define $s_1(a)=\# E^1(a)$.

The block of exceptional hypersurfaces $E^1(a)$ intervenes
in our desingularization algorithm in a way that can be thought
of intuitively as follows.
(The idea will be made precise in Sections 2 and 3.)
The exceptional hypersurfaces passing through $a$ but not in
$E^1(a)$ have accumulated during the recent part of our history,
when the order $\nu_1$ has not changed;
we have good control over these hypersurfaces.
But those in $E^1(a)$ accumulated long ago; we have forgotten
a lot about them in the form of our equations (for example,
if we restrict the equations of $X$ to these hypersurfaces,
their orders might increase) and we recall
them using $s_1(a)$.

In general, consider a sequence of blowings-up (1.8) with
$(r+\frac{1}{2})$-admissible centres.
($\inv_{r+\frac{1}{2}}$ can only decrease over such a sequence;
see Section 3 and Theorem B.)
Suppose that $i$ is the smallest index $k$ such that 
$\inv_{r+\frac{1}{2}}(a)= \inv_{r+\frac{1}{2}}(a_k)$.
Let $E^{r+1}(a)=\{ \, H\in E(a)\bs \bigcup_{q\le r} E^q(a):$
$H$ is transformed from $E(a_i)\,\}$.
We define $s_{r+1}(a) = \# E^{r+1}(a)$.

It is less straightforward to define the multiplicities $\nu_2(a),
\nu_3(a),\ldots$ and to show they are invariants.
Our definition depends on a construction in local coordinates
that we present in Section 3.
But we first try to convey the idea by describing the origin
of our algorithm.

\section{The origin of our approach}

Consider a hypersurface $X$, defined locally by an equation
$f(x)=0$.
Let $a\in X$ and let $d=d(a)$ denote the order of $X$ (or of $f$)
at $a$;
i.e., $d=\nu_1(a)=\mu_a(f)$.
We can choose local coordinates $(x_1,\ldots,x_n)$ in which
$a=0$ and $(\partial^d f/\partial x_n^d)(a)\neq 0$; then we can
write
\[
f(x) = c_0(\tx) + c_1(\tx) x_n + \cdots + c_{d-1} (\tx) x_n^{d-1}
+c_d (x) x_n^d
\]
in a neighbourhood of $a$, where $c_d(x)$ does not vanish.
($\tx$ means $(x_1,\ldots,x_{n-1})$.)
Assume for simplicity that $c_d(x)\equiv 1$ (for example,
by the Weierstrass preparation theorem, but see Remark 2.3 below).
We can also assume that $c_{d-1}(\tx)\equiv 0$, by ``completing
the $d$'th power'' (i.e., by the coordinate change $x_n'=x_n+c_{d-1}
(\tx)/d$); thus
$$
f(x) = c_0(\tx) + \cdots + c_{d-2}(\tx) x_n^{d-2} + x_n^d .
\leqno(2.1)
$$

Our aim is to simplify $f$ by blowing up with smooth centre
in the {\em equimultiple locus\/} of $a=0$; i.e., in the locus
of points of order $d$,
\[
S_{(f,d)} = \{\, x:\ \mu_x(f) = d\,\} .
\]
The representation (2.1) makes it clear that the equimultiple
locus lies in a smooth subspace of codimension $1$; in fact,
by elementary calculus,
$$
S_{(f,d)} = \{\, x:\ x_n=0\ \mbox{and}\ \mu_\tx (c_q)\ge d-q, \ 
q = 0,\ldots,d-2\,\} .
\leqno(2.2)
$$
The idea now is that the given data $\big( f(x),d\big)$ involving
$n$ variables should be equivalent, in some sense,
to the data $\cH_1(a)=\big\{ \big( c_q(\tx),d-q\big)\big\}$ in
$n-1$ variables,
thus making possible an induction on the number of variables.
(Here in ``year zero'', before we begin to blow up,
$\nu_2(a) = \min_q \mu_a (c_q) / (d-q)$.)

\rem {\em 2.3.}\quad 
For the global desingularization algorithm, the Weierstrass
preparation theorem must be avoided for two important reasons:
(1) It may take us outside the given category (for example,
in the algebraic case).
(2) Even in the complex-analytic case, we need to prove that
$\inv_X$ is semicontinuous in the sense that any point admits
a coordinate neighbourhood $V$ such that, given $a\in V$,
$\{\, x\in V:\ \inv_X(x)\le \inv_X(a)\,\}$ is Zariski-open
in $V$ (i.e., is the complement of a closed analytic subset).
We therefore need a representation like (2.2) that is valid
in a Zariski-open neighbourhood of $a$ in $V$.
This can be achieved in the following simple way that
involves neither making $c_d(x)\equiv 1$ nor explicitly
completing the $d$'th power:
By a linear coordinate change, we can assume that
$(\partial^d f/\partial x_n^d)(a)\ne 0$.
Then in the Zariski-open neighbourhood of $a$ where
$(\partial^d f/\partial x_n^d)(x)\ne 0$, we let $N_1=N_1(a)$
denote the submanifold of codimension one (in our category)
defined by $z=0$, where $z=\partial^{d-1} f/\partial x_n^{d-1}$,
and we take $\cH_1(a)=\big\{ \,\big( (\partial^q f/\partial x_n^q)
|N_1,\ d-q\big)\,\big\}$.
As before, we have $S_{(f,d)}=\{\, x:\ x\in N_1$ and 
$\mu_x(h)\ge \mu_h$, for all 
$(h,\mu_h)=\big( (\partial^q f/\partial x_n^q)|N_1,\ d-q\big)
\in \cH_1(a)\,\}$.
\medskip

We now consider the effect of a blowing-up $\sigma$ with smooth
centre $C\subset S_{(f,d)}$.
By a transformation of the variables $(x_1,\ldots,x_{n-1})$,
we can assume that in our local coordinate neighbourhood $U$
of $a$, $C$ has the form
$$
Z_I = \{\, x:\ x_n=0\ \hbox{ and }\ x_i=0,\ i\in I\,\} ,
\leqno(2.4)
$$
where $I\subset \{ 1,\ldots,n-1\}$.
According to 1.2 above, $U'=\sigma^{-1} (U)$ is covered by coordinate
charts $U_i'$, $i\in I\cup \{ n\}$, where each $U_i'$ has coordinates
$y=(y_1,\ldots,y_n)$ in which $\sigma$ is given by
$$
\begin{array}{ll}
x_i = y_i      &\\
x_j = y_i y_j ,&\qquad j\in (I\cup\{ n\})\bs \{i\}, \\
x_j = y_j\ \, ,&\qquad j\not\in I\cup \{ n\}.
\end{array}
$$
In each $U_i'$, we can write $f\big(\sigma(y)\big) = y_i^d f'(y)$;
the strict transform $X'$ of $X$ by $\sigma$ is defined in
$U_i'$ by the equation $f'(y)=0$.
(To be as simple as possible, we continue to assume $c_d(x)\equiv 1$,
though we could just as well work with the set-up of Remark 2.3;
see [BM5, Prop. 4.12].)
By (2.1), if $i\in I$, then
$$
f'(y) = c_0'(\ty) + \cdots + c_{d-2}' (\ty) y_n^{d-2} + y_n^d ,
\leqno(2.5)
$$
where
$$
c_q'(\ty) = y_i^{-(d-q)} c_q\big(\tsigma(\ty)\big),\qquad
q=0,\ldots,d-2.
\leqno(2.6)
$$
The analogous formula for the strict transform in the chart
$U_n'$ shows that $f'$ is invertible at every point of $U_n'\bs
\bigcup_{i\in I} U_i' = \{ \, y\in U_n':\ y_i=0,\ i\in I\, \}$;
in other words, $X'\cap U'\subset \bigcup_{i\in I} U_i'$.

The formula for $f'(y)$ above shows that the representation
(2.2) of the equimultiple locus (or that of Remark 2.3) is
stable under $\nu_1$-admissible blowing up when the order
does not decrease; i.e., at a point $a'\in U_i'$ where
$d(a')=d$, $S_{(f',d)} = \{ \, y:\ y_n=0$ and $\mu_\ty(c_q')\ge
d-q$, $q=0,\ldots,d-2\,\}$, where $N_1(a')=\{y_n=0\}$ is the strict
transform of $N_1(a)=\{ x_n=0\}$ and the $c_q'$ are given by
the transformation law (2.6).
The latter is not strict transform, but something intermediate
between strict and total transform $c_q\circ\sigma$.
It is essentially for this reason that some form of embedded
desingularization will be needed for the coefficients $c_q$
(i.e., in the inductive step) even to prove a weaker form of
resolution of singularities for $f$.

$N_1(a)$ is called a smooth hypersurface of {\em maximal contact}
with $X$; this means a smooth hypersurface that contains the
equimultiple locus of $a$, stably (i.e., even after admissible
blowings-up as above).
The existence of $N_1(a)$ depends on characteristic zero.
A maximal  contact hypersurface is crucial to our construction
by increasing codimension.
(In 1.12 above, $E^1(a)$ is the block of exceptional hypersurfaces
that do not necessarily have normal crossings with respect to
a maximal contact hypersurface; the term $s_1(a)$ in $\inv_X(a)$
is needed to deal with these exceptional divisors.)

We will now make a simplifying assumption on the
coefficients $c_q$:
Let us assume that one of these functions is a monomial
(times an invertible factor) that divides all the others,
but in a way that respects the different ``multiplicities''
$d-q$ associated with the transformation law (2.6);
in other words, let us make the monomial assumption on the
$c_q^{1/(d-q)}$ (to equalize the ``assigned multiplicities'' $d-q$)
or on the $c_q^{d!/(d-q)}$ (to avoid fractional powers).
We assume, then, that
$$
c_q(\tx)^{d!/(d-q)} = (\tx^\Omega)^{d!} c_q^* (\tx),\qquad
q=0,\ldots,d-2,
\leqno(2.7)
$$
where $\Omega=(\Omega_1,\ldots,\Omega_{n-1})$ with $d!\Omega_i\in\BN$
for each $i$, 
$\tx^\Omega = x_1^{\Omega_1}\cdots x_{n-1}^{\Omega_{n-1}}$,
and the $c_q^*$ are regular functions on $\{ x_n=0\}$ such that
$c_q^*(a)\neq 0$ for some $q$.
We also write $\Omega=\Omega(a)$.

We can regard (2.7) provisionally as an assumption made to see what
happens in a simple test case, but in fact we can reduce to this case
by a suitable induction on dimension (as we will see below).
(Assuming (2.7) in year zero, $\nu_2(a)=|\Omega|$, where
$|\Omega|=\Omega_1+\cdots+\Omega_{n-1}$.
But from the viewpoint of our algorithm for canonical desingularization
as presented in Section 3, the argument following is analogous
to a situation where the variables $x_i$ occurring in $\tx^\Omega$
are exceptional divisors in $E(a)\backslash E^1(a)$; in this
context, $|\Omega|$ is an invariant we call $\mu_2(a)$
(Definition 3.2) and $\nu_2 (a)=0$.)

Now, by (2.2) and (2.7),
\[
S_{(f,d)} = \{\, x:\ x_n=0\ \mbox{and}\ \mu_\tx(\tx^\Omega)\ge 1\,\} .
\]
(The order of a monomial with rational exponents has the obvious
meaning.)
Therefore (using the notation (2.4)), $S_{(f,d)} = \bigcup Z_I$,
where $I$ runs over the {\em minimal\/} subsets of $\{1,\ldots,n-1\}$
such that $\sum_{j\in I} \Omega_j \ge 1$; i.e., where $I$ runs over
the subsets of $\{ 1,\ldots,n-1\}$ such that
$$
0 \le \sum_{j\in I} \Omega_j -1 < \Omega_i,\qquad 
\mbox{ for all $i\in I$}.
\leqno(2.8)
$$

Consider the blowing-up $\sigma$ with centre $C=Z_I$, for one such
$I$.
By (2.7), in the chart $U_i'$ we have
$$
c_q'(\ty)^{d!/(d-q)} = \big( y_1^{\Omega_1}\cdots
y_i^{\sum_I\Omega_j-1}\cdots y_{n-1}^{\Omega_{n-1}}\big)^{d!}
c_q^* \big(\tsigma(\ty)\big) ,
\leqno(2.9)
$$
$q=0,\ldots,d-2$.
Suppose $a'\in \sigma^{-1} (a)\cap U_i'$.
By (2.5), $d(a')\le d(a)$.
Moreover, if $d(a')=d(a)$, then by (2.8) and (2.9), $1\le |\Omega(a')|
< |\Omega(a)|$.
In particular, the order $d$ must decrease after at most $d!|\Omega|$
such blowings-up.

The question then is whether we can reduce to the hypothesis (2.7)
by induction on dimension, replacing $(f,d)$ in some sense by
the collection $\cH_1(a)=\{(c_q,d-q)\}$ on the submanifold
$N_1=\{x_n=0\}$.
To set up the induction, we would have to treat from the start
a collection $\cF_1=\{ (f,\mu_f)\}$ rather than a single pair $(f,d)$.
(A general $X$ is, in any case, defined locally by several equations.)
Moreover, since the transformation law (2.6) is not strict transform,
we would have to reformulate the original problem to not only
desingularize $X$: $f(x)=0$, but also make its total transform
normal crossings.
To this end, suppose that $f(x)=0$ actually represents the strict
transform of our original hypersurface in that year in the history
of the blowings-up involved where the order at $a$ first becomes $d$.
(We are following the transforms of the hypersurface at a sequence
of points ``$a$'' over some original point.)
Suppose there are $s=s(a)$ accumulated exceptional hypersurfaces
$H_p$ passing through $a$; as above, we can also assume that $H_p$
is defined near $a$ by an equation
\[
x_n + b_p(\tx) = 0,
\]
$1\le p\le s$.
(Each $\mu_a(b_p)\ge 1$.)
The transformation law for the $b_p$ analogous to (2.6) is
\[
b'_p(\ty) = y_i^{-1} b_p \big(\tsigma (\ty)\big),\qquad
p=1,\ldots,s.
\]
Suppose now that in (2.7) we also have
\[
b_p(\tx)^{d!} = (\tx^\Omega)^{d!} b_p^* (\tx),\qquad
p=1,\ldots,s
\]
(and assume that either some $c_q^*(a)\neq 0$ or some
$b_p^*(a)\neq 0$).
Then the argument above shows that $\big( d(a'),s(a')\big)\le
\big( d(a),s(a)\big)$ (with respect to the lexicographic
ordering of pairs), and that if $\big( d(a'),s(a')\big)=
\big( d(a),s(a)\big)$ then $1\le |\Omega(a')| < |\Omega(a)|$.
($s(a')$ counts the exceptional hypersurfaces $H_p'$ passing
through $a'$.
As long as $d$ does not drop, the new exceptional hypersurfaces
accumulate simply as $y_i=0$ for certain $i=1,\ldots,n-1$,
in suitable coordinates $(y_1,\ldots,y_{n-1})$ for the strict
transform $N'=\{y_n=0\}$ of $N=\{ x_n=0\}$.)

The induction on dimension can be realized in various ways.
The simplest --- the method of [BM1, Section 4] --- is to apply
the inductive hypothesis within a coordinate chart to the function
of $n-1$ variables given by the product of all nonzero
$c_q^{d!/(d-q)}$, all nonzero $b_p^{d!}$, and all their nonzero
differences.
The result is (2.7) and (2.10) (with $c_q^*(a)\neq 0$ or
$b_p^*(a)\neq 0$ for some $q$ or $p$; see [BM1, Lemma 4.7]).
Pullback of the coefficients $c_q$ by a blowing-up in
$(x_1,\ldots,x_{n-1})$ with smooth centre $C$, corresponds to
strict transform of $f$ by the blowing-up with centre
$C\times \{ x_n-\mbox{axis}\}$.
Thus we sacrifice the condition that each centre lie in the 
equimultiple locus (or even in $X$!).
But we do get a very simple proof of local uniformization.
In fact, we get the conclusion (2) of our desingularization
theorem A, using a mapping $\sigma$: $M'\rightarrow M$ which is a
composite of mappings that are either blowings-up with smooth
centres or surjections of the form $\coprod_j U_j \rightarrow\bigcup_j
U_j$, where the latter is a locally-finite open covering of a manifold
and $\coprod$ means disjoint union.

To prove our canonical desingularization theorem, we repeat the
construction above in increasing codimension to obtain
$\inv_X(a)=\big(\nu_1(a)$, $s_1(a);\, \nu_2(a),\ldots\,\big)$ ---
$\big(\nu_1(a),s_1(a)\big)$ is $\big( d(a),s(a)\big)$ above --- 
together with a corresponding local ``presentation''.
The latter means a local description of the locus of constant values
of the invariant in terms of regular functions with assigned
multiplicities, that survives certain blowings-up.
($N_1(a), \cH_1(a)$ above is a presentation of $\nu_1$ at $a$.)

\section{The desingularization algorithm}

In this section we give a constructive definition of $\inv_X$
together with a corresponding presentation (in the hypersurface
case).
We illustrate the construction by applying the desingularization
algorithm to an example --- a surface whose desingularization
involves all the features of the general hypersurface case.
We will use horizontal lines to separate from the example
the general considerations that are needed at each step.

\ex {\em 3.1.}\quad
Let $X\subset \BK^3$ denote the hypersurface $g(x)=0$, where
$g(x)=x_3^2 -x_1^2 x_2^3$.

\begin{center}
{\hskip .2in{\epsfxsize=2.5in\epsfbox{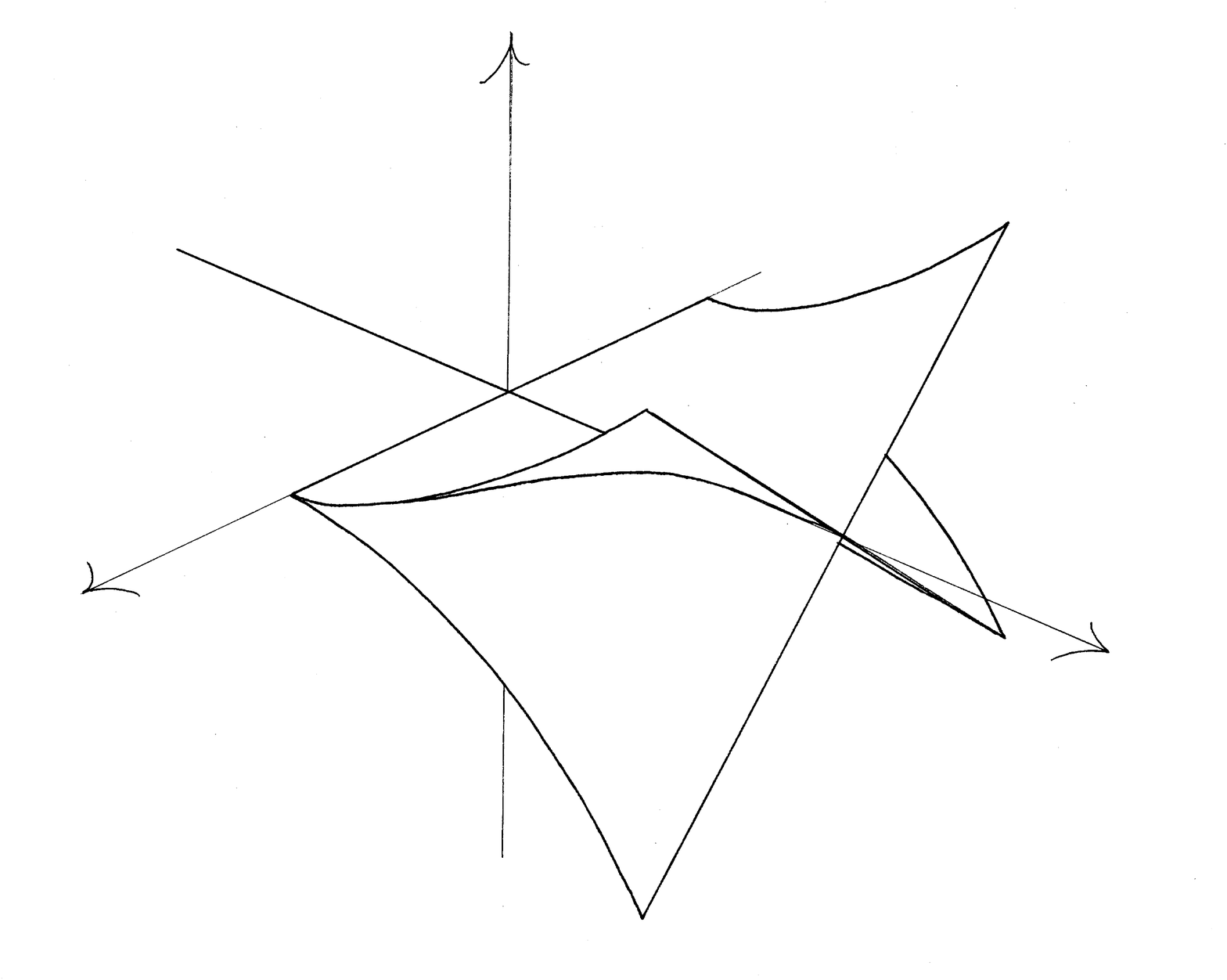}{\vskip -2.25in\hskip -.11in{$x_3$}}
{\vskip .35in\hskip 2in{$X$}}
{\vskip .07in\hskip -.88in{$a=0$}}
{\vskip .25in\hskip -2.04in{$x_1$}}
{\vskip -.13in\hskip 2.47in{$x_2$}}}}
\end{center}
\vspace{.75in}

Let $a=0$.
Then $\nu_1(a)=\mu_a(g)=2$.
Of course, $E(a)=\emptyset$, so that $s_1(a)=0$.
(This is ``year zero''; there are no exceptional hypersurfaces.)
Thus $\inv_1(a)=\big(\nu_1(a),s_1(a)\big)=(2,0)$.
Let $\cG_1(a)=\{ (x_3^2-x_1^2 x_2^3,2)\}$.
We say that $\cG_1(a)$ is a {\em codimension\/} 0 {\em presentation
of\/} $\inv_{\frac{1}{2}}=\nu_1$ {\em at\/} $a$.
(Here where $s_1(a)=0$, we can also say that $\cG_1(a)$ is a
codimension $0$ presentation of $\inv_1=(\nu_1,s_1)$ at $a$.)

\linee

In general, consider a hypersurface $X\subset M$.
Let $a\in M$ and let $S_{\inv_{\frac{1}{2}}}(a)$ denote the germ at 
$a$ of $\{\, x:\ \inv_{\frac{1}{2}}(x)\ge\inv_{\frac{1}{2}}(a)\,\}$ 
$=$ the germ at $a$ of $\{\, x:\ \inv_{\frac{1}{2}}(x)=
\inv_{\frac{1}{2}}(a)\,\}$.
If $g\in\cO_{M,a}$ generates the local ideal $\cI_{X,a}$ of $X$
and $d=\nu_1(a)=\mu_a(g)$, then $\cG_1(a)=\{(g,d)\}$ is a
codimension $0$ presentation of $\inv_{\frac{1}{2}}=\nu_1$ at $a$.
This means $S_{\inv_{\frac{1}{2}}}(a)$ coincides with the germ of
the ``equimultiple locus'' $S_{\cG_1(a)} = \{\, x:\ \mu_x(g)=d\,\}$,
and that the latter condition survives certain transformations.

More generally, suppose that $\cG_1(a)$ is a finite collection
of pairs $\{ (g,\mu_g)\}$, where each $g$ is a germ at $a$ of
a regular function (i.e., $g\in\cO_{M,a}$) with an ``assigned
multiplicity'' $\mu_g\in\BQ$, and where we assume that 
$\mu_a(g)\ge\mu_g$ for every $g$.
Set
\[
S_{\cG_1(a)}\ =\ \{\, x:\ \mu_x(g)\ge\mu_g,\ \mbox{for all}\ 
(g,\mu_g)\in\cG_1(a)\,\} ;
\]
$S_{\cG_1(a)}$ is well-defined as a germ at $a$.
To say that $\cG_1(a)$ is a {\em codimension\/} $0$ 
{\em presentation of\/} $\inv_{\frac{1}{2}}$ 
{\em at\/} $a$ means that
\[
S_{\inv_{\frac{1}{2}}}(a) = S_{\cG_1(a)}
\]
and that this condition survives certain transformations:

To be precise, we will consider triples of the form $\big( N=N(a)$,
$\cH(a)$, $\cE(a)\big)$, where:

$N$ is a germ of a submanifold of codimension $p$ at $a$
(for some $p\ge 0$).

$\cH(a)=\{(h,\mu_h)\}$ is a finite collection of pairs $(h,\mu_h)$,
where $h\in\cO_{N,a}$, $\mu_h\in\BQ$ and $\mu_a(h)\ge\mu_h$.

$\cE(a)$ is a finite set of smooth (exceptional) hyperplanes
containing $a$, such that $N$ and $\cE(a)$ simultaneously have
normal crossings and $N\not\subset H$, for all $H\in \cE(a)$.

A {\em local blowing-up\/} $\sigma$: $M'\rightarrow M$ over a 
neighbourhood $W$ of $a$, with smooth centre $C$, means the 
composite of a blowing-up $M'\rightarrow W$ with centre $C$, and 
the inclusion $W\hookrightarrow M$.

\defn {\em 3.2.}\quad We say that $\big( N(a),\cH(a),\cE(a)\big)$
is a {\em codimension\/} $p$ {\em presentation of\/} 
$\inv_{\frac{1}{2}}$ {\em at\/} $a$ if:

\smallskip

(1) $S_{\inv_{\frac{1}{2}}}(a) = S_{\cH(a)}$, where 
$S_{\cH(a)}=\{\, x\in N:\ \mu_x(h)\ge\mu_h$, for all
$(h,\mu_h)\in\cH(a)\,\}$ (as a germ at $a$).

(2) Suppose that $\sigma$ is a $\frac{1}{2}$-admissible local 
blowing-up at $a$ (with smooth centre $C$).
Let $a'\in\sigma^{-1}(a)$.
Then $\inv_{\frac{1}{2}}(a')=\inv_{\frac{1}{2}}(a)$ if and only if 
$a'\in N'$ (where $N'=N(a')$ denotes the strict transform of $N$)
and $\mu_{a'}(h')\ge \mu_{h'}$ for all $(h,\mu_h)\in\cH(a)$,
where $h'=y_\exc^{-\mu_h} h\circ\sigma$ and $\mu_{h'}=\mu_h$.
($y_\exc$ denotes a local generator of $\cI_{\sigma^{-1}(C)}$.)
In this case, we will write $\cH(a')=\{\, (h',\mu_{h'}):\ (h,\mu_h)\in
\cH(a)\,\}$ and $\cE(a')=\{\, H':\ H\in\cE(a)\,\}\cup \{\sigma^{-1}
(C)\}$.

(3) Conditions (1) and (2) continue to hold for the transforms
$X'$ and $\big( N(a'),\cH(a'),\cE(a')\big)$ of our data by sequences
of morphisms of the following three types, at points $a'$ in the
fibre of $a$ (to be also specified).

\medskip

The three types of morphisms allowed are the following.
(Types (ii) and (iii) are not used in the actual desingularization
algorithm.
They are needed to prove invariance of the terms $\nu_2(a), \nu_3(a),
\ldots$ of $\inv_X(a)$ by making certain sequences of ``test
blowings-up'', as we will explain in Section 4; they are not
explicitly needed in this section.)
\medskip

(i) $\frac{1}{2}$-{\em admissible local blowing-up\/} $\sigma$, and
$a'\in\sigma^{-1}(a)$ such that 
$\inv_{\frac{1}{2}}(a')=\inv_{\frac{1}{2}}(a)$.

(ii) {\em Product with a line.\/}  $\sigma$ is a projection
$M'=W\times\BK\rightarrow W\hookrightarrow M$, where $W$ is a
neighbourhood of $a$, and $a'=(a,0)$.

(iii) {\em Exceptional blowing-up.\/}
$\sigma$ is a local blowing-up $M'\rightarrow W\hookrightarrow M$ over
a neighbourhood $W$ of $a$, with centre $H_0\cap H_1$, where
$H_0,H_1\in\cE(a)$, and $a'$ is the unique point of
$\sigma^{-1}(a)\cap H_1'$.

\medskip

The data is transformed to $a'$ in each case above, as follows:

\medskip

(i) $X'=$ strict transform of $X$;
$\big( N(a'),\cH(a'),\cE(a')\big)$ as defined in 3.2(2) above.

(ii) and (iii) $X'=\sigma^{-1}(X)$, $N(a')=\sigma^{-1}(N)$,
$\cH(a')=\{(h\circ\sigma,\mu_h)\}$.
$\cE(a') = \{\,\sigma^{-1}(H):\ H\in\cE(a)\,\}\cup \{W\times 0\}$
in case (ii); $\cE(a')=\{\, H':\ H\in\cE(a),\ a'\in H'\,\}\cup
\{\sigma^{-1}(C)\}$ in case (iii).

\medskip

If $\big( N(a),\cH(a),\cE(a)\big)$ is a presentation of
$\inv_{\frac{1}{2}}$ at $a$, then $N(a)$ is called a subspace
of {\em maximal contact} (cf. Section 2).

Suppose now that $\cG_1(a)$ is a codimension $0$ presentation
of $\inv_{\frac{1}{2}}$ at $a$.
(Implicitly, $N(a)=M$ and $\cE(a)=\emptyset$.)
Assume, moreover, that there exists $(g,\mu_g)=(g_*,\mu_{g_*})\in
\cG_1(a)$ such that $\mu_a(g_*)=\mu_{g_*}$ (as in Example 3.1).

We can always assume that each $\mu_g\in\BN$, and even that all
$\mu_g$ coincide:
Simply replace each $(g,\mu_g)$ by $(g^{e/\mu_g},e)$, for
suitable $e\in\BN$.

Then, after a linear coordinate change if necessary, we can assume
that $(\partial^d g_*/\partial x_n^d)(a)\neq 0$, where $d=\mu_{g_*}$.
Set
\begin{eqnarray*}
z & = & \frac{\partial^{d-1} g_*}{ \partial x_n^{d-1}}\in \cO_{M,a}\\
N_1 = N_1(a) & = &\{ z=0\} \\
\cH_1(a) & = &\left\{ \left( \frac{\partial^q g}{\partial x_n^q} \bigg|
_{N_1}, \mu_g-q\right):\ 0\le q < \mu_g,\ (g,\mu_g)\in\cG_1(a)\right\} .
\end{eqnarray*}
Then $\big( N_1(a),\cH_1(a),\cE_1(a)=\emptyset\big)$ is a codimension
$1$ presentation of $\inv_{\frac{1}{2}}$ at $a$:
This is an assertion about the way our data transforms under
sequences of morphisms of types (i), (ii) and (iii) above.
The effect of a transformation of type (i) is essentially described
by the calculation in Section 2.
The effect of a transformation of type (ii) is trivial, and
that for type (iii) can be understood in a similar way to (i):
see [BM5, Props. 4.12 and 4.19] for details.

\defn {\em 3.3.}\quad
We define
\[
\mu_2(a) = \min_{\cH_1(a)}\, \frac{\mu_a(h)}{\mu_h}.
\]

Then $1\le \mu_2(a)\le\infty$.
If $\cE(a)=\emptyset$ (as in year zero), we set
\[
\nu_2(a)=\mu_2(a)
\]
and $\inv_{1\frac{1}{2}}(a)=\big(\inv_1(a);\nu_2(a)\big)$.
Then $\nu_2(a)\le\infty$.
Moreover, $\nu_2(a)=\infty$ if and only if $\cG_1(a)\sim \{(z,1)\}$.
(This means that the latter is also a presentation of 
$\inv_{\frac{1}{2}}$ at $a$.)
If $\nu_2(a)=\infty$, then we set $\inv_X(a)=\inv_{1\frac{1}{2}}(a)$.
$\inv_X(a)=(d,0,\infty)$ if and only if $X$ is defined (near $a$)
by the equation $z^d=0$; in this case, the desingularization
algorithm can do no more, unless we blow-up with centre
$|X|$!

\linee

In Example 3.1, $\mu_a(g)=2=\mu_g$, and by the construction
above we get the following codimension $1$ presentation
of $\inv_{\frac{1}{2}}$ (or $\inv_1$) at $a$:
\[
N_1(a)=\{x_3=0\},\qquad
\cH_1(a) = \{ (x_1^2 x_2^3,2)\}.
\]
Thus $\nu_2(a)=\mu_2(a)=5/2$.
As a codimension $1$ presentation of $\inv_{1\frac{1}{2}}$
(or $\inv_2$) at $a$, we can take
\[
N_1(a),\qquad
\cG_2(a) = \{ (x_1^2 x_2^3,5)\} .
\]

\linee

In general, ``presentation of $\inv_r$'' (or ``of
$\inv_{r+\frac{1}{2}}$'') means the analogue of ``presentation
of $\inv_{\frac{1}{2}}$'' above.
Suppose that $\big( N_1(a),\cH_1(a)\big)$ is a codimension $1$
presentation of $\inv_1$ at $a$ $\big(\cE_1(a)=\emptyset\big)$.
Assume that $1\le\nu_2(a)<\infty$.
(In year zero, we always have $\nu_2(a)=\mu_2(a)\ge 1$.)
Let
\[
\cG_2(a) = \big\{\, \big(h,\nu_2(a)\mu_h\big):\ 
(h,\mu_h)\in\cH_1(a)\big\} .
\]
Then $\big( N_1(a),\cG_2(a)\big)$ is a codimension $1$ presentation
of $\inv_{\frac{1}{2}}$ at $a$ (or of $\inv_2$ at $a$, when
$s_2(a)=0$ as here).
Clearly, there exists $(g_*,\mu_{g_*})\in \cG_2(a)$ such that
$\mu_a(g_*)=\mu_{g_*}$.

This completes a cycle in the recursive
definition of $\inv_X$, and we can now repeat the above constructions:
Let $d=\mu_{g_*}$.
After a linear transformation of the coordinates 
$(x_1,\ldots,x_{n-1})$ of $N_1(a)$, we can assume that 
$(\partial^d g_* /\partial x_{n-1}^d) (a)\neq 0$.
We get a codimension $2$ presentation of $\inv_2$ at $a$ by taking
\begin{eqnarray*}
N_2(a) & = & \left\{\, x\in N_1(a):\ \frac{\partial^{d-1} g_*}
{\partial x_{n-1}^{d-1}} (x) = 0\,\right\} ,\\
\cH_2(a) & = & \left\{\,\left( \frac{\partial^q g}{\partial x_{n-1}^q}
\bigg|_{N_2(a)} , \mu_g -q\right):\ 0\le q<\mu_g,\ (g,\mu_g)\in
\cG_2(a)\,\right\} .
\end{eqnarray*}

In our example, the calculation of a codimension $2$ presentation
can be simplified by the following useful observation:
Suppose there is $(g,\mu_g)\in\cG_2(a)$ with $\mu_a(g)=\mu_g$
and $g=\prod g_i^{m_i}$.
If we replace $(g,\mu_g)$ in $\cG_2(a)$ by the collection of
$(g_i,\mu_{g_i})$, where each $\mu_{g_i}=\mu_a(g_i)$,
then we obtain an (equivalent) presentation of $\inv_2$.

\linee

In our example, therefore,
\[
N_1(a) = \{ x_3=0\},\qquad
\cG_2(a) = \{ (x_1,1), (x_2,1)\}
\]
is a codimension $1$ presentation of $\inv_2$ at $a$.
It follows immediately that
\[
N_2(a)=\{x_2=x_3=0\},\qquad
\cH_2(a) = \{ (x_1,1)\}
\]
is a codimension $2$ presentation of $\inv_2$ at $a$.
Then $\nu_3(a)=\mu_3(a)=1$ and, as a codimension $3$
presentation of $\inv_{2\frac{1}{2}}$ (or of $\inv_3$) at $a$,
we can take
\[
N_3(a) = \{ x_1=x_2=x_3=0\},\qquad
\cH_3(a) = \emptyset .
\]
We put $\nu_4(a)=\mu_4(a)=\infty$.
Thus we have
\[
\inv_X (a) = (2,0;\, 5/2,0;\, 1,0;\,\infty)
\]
and $S_{\inv_X}(a)=S_{\inv_3}(a)=N_3(a)=\{a\}$.
The latter is the centre $C_0$ of our first blowing-up
$\sigma_1$: $M_1\rightarrow M_0=\BK^3$; $M_1$ can be covered
by three coordinate charts $U_i$, $i=1,2,3$, where each
$U_i$ is the complement in $M_1$ of the strict transform
of the hyperplane $\{ x_i=0\}$.
The strict transform $X_1=X'$ of $X$ lies in $U_1\cup U_2$.
To illustrate the algorithm, we will follow our construction
at a sequence of points over $a$, choosing after each
blowing-up a point in the fibre where $\inv_X$ has a
maximum value in a given coordinate chart.

\medskip

{\em Year one.}\quad 
$U_1$ has a coordinate system $(y_1,y_2,y_3)$ in which 
$\sigma_1$ is given by the transformation
\[
x_1=y_1,\quad
x_2=y_1y_2,\quad
x_3=y_1y_3 .
\]
Then $X_1\cap U_1 =V(g_1)$, where
\[
g_1 = y_1^{-2} g\circ\sigma_1 = y_3^2 - y_1^3 y_2^3 .
\]
Consider $b=0$.
Then $E(b)=\{ H_1\}$, where $H_1$ is the exceptional
hypersurface $H_1=\sigma_1^{-1}(a) = \{ y_1=0\}$.
Now, $\nu_1(b)=2=\nu_1(a)$.
Therefore $E^1(b)=\emptyset$ and $s_1(b)=0$.
We write $\cE_1(b)=E(b)\bs E^1(b)$, so that $\cE_1(b)=E(b)$
here.
Let $\cF_1(b)=\cG_1(b)=\{ (g_1,2)\}$.
Then $\big( N_0(b)=M_1, \cF_1(b),\cE_1(b)\big)$ is a
codimension $0$ presentation of $\inv_1$ at $b$.
Set
\[
N_1(b) = \{ y_3=0\} = N_1(a)',\qquad
\cH_1(b) = \{ (y_1^3 y_2^3,2)\} ;
\]
$\big( N_1(b), \cH_1(b),\cE_1(b)\big)$ is a codimension $1$
presentation of $\inv_1$ at $b$.
As before,
\[
\mu_2(b) = \min_{\cH_1(b)} \frac{\mu_b(h)}{\mu_h} =
\frac{6}{2} = 3 .
\]
But, here, in the presence of nontrivial $\cE_1(b)$,
$\nu_2(b)$ will involve first factoring from the $h\in
\cH_1(b)$ the exceptional divisors in $\cE_1(b)$
(taking, in a sense, ``internal strict transforms'' at
$b$ of the elements of $\cH_1(a)$).

\linee

In general, we define
\[
\cF_1(b) = \cG_1(b)\cup \big(E^1(b),1\big) ,
\]
where $\big( E^1(b),1\big)$ denotes $\{ \,(y_H,1):\ H\in E^1(b)\,\}$,
and $y_H$ means a local generator of the ideal of $H$.
Then $\big( N_0(b),\cF_1(b),\cE_1(b)\big)$ is a codimension $0$
presentation of $\inv_1=(\nu_1,s_1)$ at $b$, and there is a
codimension $1$ presentation $\big( N_1(b),\cH_1(b),\cE_1(b)\big)$ as
before.
The construction of Section 2 above shows that we can choose
the coordinates $(y_1,\ldots,y_{n-1})$ of $N_1(b)$ so that
each $H\in\cE_1(b)=E(b)\bs E^1(b)$ is $\{y_i=0\}$, for some $i=
1,\ldots,n-1$; we again write $y_H=y_i$.
(In other words, $\cE_1(b)$ and $N_1(b)$ simultaneously have
normal crossings, and $N_1(b)\not\subset H$, for all $H\in\cE_1(b)$.)

\defn {\em 3.4.}\quad 
For each $H\in \cE_1(b)$, we set
\[
\mu_{2H}(b) = \min_{(h,\mu_h)\in\cH_1(b)} \frac
{\mu_{H,b}(h)}{\mu_h} ,
\]
where $\mu_{H,b}(h)$ denotes the {\em order\/} of $h$ {\em along\/}
$H$ at $b$;
i.e., the order to which $y_H$ factors from $h\in\cO_{N,b}$,
$N=N_1(b)$, or $\max\{\, k: \  h\in\cI_{H,b}^k\,\}$, where
$\cI_{H,b}$ is the ideal of $H\cap N$ in $\cO_{N,b}$.
We define
\[
\nu_2(b) = \mu_2(b) - \sum_{H\in\cE_1(b)} \mu_{2H} (b) .
\]

\linee

In our example,
\[
\nu_2(b) = \mu_2(b) - \mu_{2H_1}(b) = 3 - \frac{3}{2} = \frac{3}{2} .
\]

\linee

Write
\[
D_2(b) = \prod_{H\in \cE_1(b)} y_H^{\mu_{2H}(b)} .
\]
Suppose, as before, that all $\mu_h$ are equal: say
all $\mu_h = d\in\BN$.
Then $D^d=D_2(b)^d$ is the greatest common divisor
of the $h$ that is a monomial in the exceptional
coordinates $y_H$, $H\in\cE_1(b)$.
For each $h\in \cH_1(b)$, write $h=D^d g$ and set
$\mu_g = d\nu_2 (b)$; then $\mu_b(g)\ge\mu_g$.
Clearly, $\nu_2(b) = \min_g \mu_b (g)/d$.
Moreover, $0\le \nu_2(b)\le\infty$, and $\nu_2(b)=\infty$
if and only if $\mu_2(b)=\infty$.

If $\nu_2(b)=0$ or $\infty$, we put 
$\inv_X(b)=\inv_{1\frac{1}{2}}(b)$.
If $\nu_2(b)=\infty$, then $S_{\inv_X}(b)=N_1(b)$.
If $\nu_2(b)=0$ and $\cG_2(b) = \big\{\big( D_2(b),1\big)\big\}$,
then $\big( N_1(b),\cG_2(b),\cE_1(b)\big)$ is a codimension
$1$ presentation of $\inv_X$ at $b$; in particular,
\[
S_{\inv_X}(b) = \big\{\, y\in N_1(b):\ \mu_y\big( D_2(b)\big)\ge
1\big\}
\]
(cf. Section 2).

Consider the case that $0<\nu_2(b)<\infty$.
Let $\cG_2(b)$ denote the collection of pairs $(g,\mu_g)=
\big( g, d\nu_2(b)\big)$ for all $(h,\mu_h)=(h,d)$, as above,
together with the pair $\big( D_2(b)^d,\big(1-\nu_2(b)\big)d\big)$
{\em provided that\/} $\nu_2(b)<1$.
Then $\big( N_1(b),\cG_2(b),\cE_1(b)\big)$ is a codimension
$1$ presentation of $\inv_{1\frac{1}{2}}$ at $b$.

In the latter case, we introduce $E^2(b)\subset \cE_1(b)$
as in 1.12, and we set $s_2(b) = \# E^2(b)$, $\cE_2(b)=
\cE_1(b)\bs E^2(b)$.
Set
\[
\cF_2(b) = \cG_2(b) \cup \big( E^2(b),1\big) .
\]
Then $\big( N_1(b),\cF_2(b),\cE_1(b)\big)$ is a codimension
$1$ presentation of $\inv_2$ at $b$, and we can pass to a
codimension $2$ presentation $\big( N_2(b),\cH_2(b),\cE_2(b)\big)$.
Here it is important to replace $\cE_1(b)$ by the subset
$\cE_2(b)$, to have the property that $\cE_2(b)$, $N_2(b)$
simultaneously have normal crossings and $N_2(b)\not\subset H$,
for all $H\in\cE_2(b)$.
(Again, the main r\^ole of $\cE$ in a presentation is
to prove invariance of the $\mu_{2H} (\cdot)$ and in general
of the $\mu_{3H}(\cdot),\ldots,$ as in Section 4.)

\linee

Returning to our example (in year one), we have
$\cH_1(b)=\{ (y_1^3 y_2^3,2)\}$, so that $D_2(b)=y_1^{3/2}$.
We can take $\cG_2(b)=\{(y_2^3,3)\}$ or, equivalently,
$\cG_2(b)=\{ (y_2,1)\}$ to get a codimension $1$ presentation
$\big( N_1(b), \cG_2(b),\cE_1(b)\big)$ of $\inv_{1\frac{1}{2}}$
at $b$.

Now, $E^2(b)=\{H_1\}$, so that $s_2(b)=1$.
We set
\[
\cF_2(b) = \cG_2(b)\cup \big( E^2(b),1\big) =
\{(y_1,1),(y_2,1)\}
\]
and $\cE_2(b)=\cE_1(b)\bs E^2(b)=\emptyset$.
Then $\big( N_1(b),\cF_2(b),\cE_1(b)\big)$ is a codimension
$1$ presentation of $\inv_2$ at $b$, and we can get
a codimension $2$ presentation $\big( N_2(b),\cH_2(b),
\cE_2(b)\big)$ of $\inv_2$ at $b$ by taking $N_2(b)=
\{ y_2=y_3=0\}$ and $\cH_2(b)=\{ (y_1,1)\}$.

It follows that $\nu_3(b)=1$.
Since $E^3(b)=\emptyset$, $s_3(b)=0$.
We get a codimension $3$ presentation of $\inv_3$ at $b$
by taking
\[
N_3(b) = \{ y_1=y_2=y_3=0\} = \{b\},\qquad
H_3(b) = \emptyset .
\]
Therefore,
\[
\inv_X(b) = \left( 2,0;\, \frac{3}{2},1;\, 1,0;\, \infty\right)
\]
and $S_{\inv_X}(b)=S_{\inv_2}(b)=\{b\}$.
The latter is the centre of the next blowing-up $\sigma_2$.
$\sigma_2^{-1}(U_1)$ is covered by 3 coordinate charts
$U_{1i}=\sigma_2^{-1}(U_1)\bs \{ y_i=0\}'$,
$i=1,2,3$.
For example, $U_{12}$ has coordinates
$(z_1,z_2,z_3)$ with respect to which $\sigma_2$
is given by
\[
y_1 = z_1 z_2,\qquad
y_2 = z_2,\qquad
y_3 = z_2 z_3 .
\]

\linee

\rem {\em 3.5.}\quad 
{\em Zariski-semicontinuity of the invariant.}
Each point of $M_j$, $j=0,1,\ldots$, admits a coordinate
neighbourhood $U$ such that, for all $x_0\in U$,
$\{ \, x\in U:\ \inv_{\!\!\lower2pt\hbox{$\displaystyle\cdot$}} (x) \le \inv_{\!\!\lower2pt\hbox{$\displaystyle\cdot$}} (x_0)\,\}$
is Zariski-open in $U$ (i.e., the complement of a Zariski-closed
subset of $U$):
For $\inv_{\frac{1}{2}}$, this is just Zariski-semicontinuity
of the order of a regular function $g$ (a local generator
of the ideal of $X$).
For $\inv_1$, the result is a consequence of the following
semicontinuity assertion for $E^1(x)$:
There is a Zariski-open neighbourhood of $x_0$ in $U$, in
which $E^1(x)=E(x)\cap E^1(x_0)$, for all $x\in S_{\inv_
{\frac{1}{2}}(x_0)} = \{\, x\in U:\ \inv_{\frac{1}{2}}(x)\ge
\inv_{\frac{1}{2}}(x_0)\,\}$.
(See [BM5, Prop. 6.6] for a simple proof.)

For $\inv_{1\frac{1}{2}}$: Suppose that $\mu_k=d\in\BN$,
for all $(h,\mu_h)\in\cH_1(x_0)$, as above.
Then, in a Zariski-open neighbourhood of $x_0$
where $S_{\inv_{1}(x)}=\{\, x:\ \inv_1(x)=\inv_1(x_0)\,\}$,
we have
\[
d\nu_2(x) = \min_{\cH_1(x_0)} \mu_x \left(
\frac{h}{D_2(x_0)^d}\right) ,\qquad
x\in S_{\inv_{1}(x_0)} .
\]
Semicontinuity of $\nu_2(x)$ is thus a consequence
of semicontinuity of the order of an element
$g=h/D_2(x_0)^d$ such that $\mu_{x_0}(g)=d\nu_2(x_0)$.

Likewise for $\inv_2$, $\inv_{2\frac{1}{2}}$, $\ldots$.

\linee

{\em Year two.}\quad 
Let $X_2$ denote the strict transform $X_1'$ of $X_1$
by $\sigma_2$.
Then $X_2\cap U_{12}=V(g_2)$, where
\[
g_2 = z_2^{-2} g_1\circ \sigma_2 = z_3^2 -z_1^3 z_2^4 .
\]
Let $c$ be the origin of $U_{12}$.
Then $E(c)=\{ H_1,H_2\}$ where
\begin{eqnarray*}
&&H_1 = \{ y_1=0\}' =\{z_1=0\} ,\\
&&H_2 = \sigma_2^{-1} (b) = \{ z_2 =0\} .
\end{eqnarray*}
We have $\nu_2(c)=2=\nu_2(a)$.
Therefore, $E^1(c)=\emptyset$, $s_1(c)=0$, $\cE_1(c)=E(c)$.
$\cF_1(c)=\cG_1(c)=\{ (g_2,2)\}$ provides a codimension $0$
presentation of $\inv_1$ at $c$, and we get a codimension $1$
presentation by taking
\[
N_1(c) = \{ z_3=0\},\qquad
\cH_1(c) = \{ (z_1^3 z_2^4,2)\} .
\]
Therefore $\mu_2(c)=7/2$, $\mu_{2H_1}(c) = 3/2$ and
$\mu_{2H_2}(c) = 4/2 =2$, so that $\nu_2(c)=0$ and
\[
\inv_X(c) = (2,0;\,0) .
\]
Moreover, $D_2(c)=z_1^{\frac{3}{2}} z_2^2$, and we get a
codimension $1$ presentation of $\inv_X=\inv_{1\frac{1}{2}}$ at $c$
using
\[
N_1(c) = \{ z_3=0\},\qquad
\cG_2(c) = \{ (z_1^{\frac{3}{2}} z_2^2,1)\} .
\]
Therefore,
\[
S_{\inv_X}(c) = S_{\inv_{1\frac{1}{2}}} (c) =
\{ z_1=z_3=0\} \cup \{ z_2=z_3=0\} ;
\]
of course, $\{ z_1=z_3=0\} = S_{\inv_X} (c)\cap H_1$
and $\{ z_2=z_3=0\} = S_{\inv_X}(c)\cap H_2$.

\linee

\rem {\em 3.6.}\quad 
In general, suppose that $\inv_X(c)=\inv_{t+\frac{1}{2}}(c)$ and
$v_{t+1}(c)=0$.
(We assume $c\in M_j$, for some $j=1,2,\ldots$.)
Then $\inv_X$ has a codimension $t$ presentation
at $c$: $N_t(c) = \{ z_{n-t+1} = \cdots = z_n = 0\}$,
$\cG_{t+1}(c) = \big\{ \big( D_{t+1} (c),1\big)\big\}$, where
$D_{t+1}(c)$ is a monomial with rational exponents in the
exceptional divisors $z_H$, $H\in \cE_t(c)$; $N_t(c)$
has coordinates $(z_1,\ldots,z_{n-t})$ in which each
such $z_H=z_i$, for some $i=1,\ldots,n-t$.
It follows that each component $Z$ of $S_{\inv_X}(c)$
has the form
\[
Z = S_{\inv_X} (c) \cap \bigcap \{\, H\in E(c):\ Z\subset H\,\} ;
\]
we will write $Z=Z_I$, where $I=\{ \,H\in E(c):\ Z\subset H\,\}$.
It follows that, if $U$ is any open neighbourhood of $c$
on which $\inv_X(c)$ is a maximum value of $\inv_X$, then
every component $Z_I$ of $S_{\inv_X}(c)$ extends to a global
smooth closed subset of $U$:

First consider any total order on $\{\, I:\ I\subset E_j\,\}$.
For any $c\in M_j$, set
\begin{eqnarray*}
J(c) & = &\max\{\, I:\ Z_I\ \mbox{is a component of}\ 
S_{\inv_X}(c)\,\} ,\\
\inv_X^\rme (c) & = & \big( \inv_X (c);\, J(c)\big) .
\end{eqnarray*}
Then $\inv_X^\rme$ is Zariski-semicontinuous (again comparing
values of $\inv_X^\rme$ lexicographically), and its
locus of maximum values on any given open subset of
$M_j$ is smooth.

Of course, given $c\in M_j$ and a component $Z_I$ of
$S_{\inv_X}(c)$, we can choose the ordering of $\{\, J:\ J\subset
E_j\,\}$ so that $I=J(c)=\max\{\, J:\ J\subset E_j\,\}$.
It follows that, if $U$ is any open neighbourhood of $c$
on which $\inv_X(c)$ is a maximum value of $\inv_X$,
then $Z_I$ extends to a smooth closed subset of $U$.

To obtain an algorithm for canonical desingularization,
we can choose as each successive centre of blowing up
the maximum locus of $\inv_X^\rme(\cdot) = \big(\inv_X(\cdot),
J(\cdot)\big)$, where $J$ is defined as above using
the following total ordering of the subsets of $E_j$:
Write $E_j=\{ H_1^j,\ldots,H_j^j\}$, where each $H_i^j$
is the strict transform in $M_j$ of the exceptional
hypersurface $H_i^i=\sigma_i^{-1} (C_{i-1})\subset
M_i$, $i=1,\ldots,j$.
We can order $\{\, I:\ I\subset E_j\,\}$ by associating
to each subset $I$ the lexicographic order of the sequence
$(\delta_1,\ldots,\delta_j)$, where $\delta_i=0$ if
$H_i^j\not\in I$ and $\delta_i=1$ if $H_i^j\in I$.

\linee

In our example, year two, we have
\[
S_{\inv_X}(c) = \big( S_{\inv_X}(c)\cap H_1\big) \cup
\big( S_{\inv_X}(c)\cap H_2\big) .
\]
(Each $H_i$ is $H_i^2$ in the notation preceding.)
The order of $H_1$ (respectively, $H_2$) is $(1,0)$
(respectively, $(0,1)$), so that $J(c)=\{ H_1\}$ and
the centre of the third blowing-up $\sigma_3$ is
$C_2=S_{\inv_X}(c)\cap H_1 = \{ z_1=z_3=0\}$.

Thus $\sigma_3^{-1}(U_{12})=U_{121}\cup U_{123}$, where
$U_{12i}=\sigma_3^{-1}(U_{12})\bs \{ z_i=0\}'$, $i=1,3$.
The strict transform of $X_2\cap U_{12}$ lies in
$U_{121}$; the latter has coordinates $(w_1,w_2,w_3)$
in which $\sigma_3$ can be written
\[
z_1=w_1,\qquad
z_2=w_2,\qquad
z_3=w_1w_3 .
\]

\medskip

{\em Year three.}\quad
Let $X_3$ denote the strict transform of $X_2$ by $\sigma_3$.
Then $X_3\cap U_{121} = V(g_3)$, where $g_3(w)=w_3^2 -
w_1 w_2^4$.
Let $d=0$ in $U_{121}$.
There are three exceptional hypersurfaces $H_1=\{ z_1=0\}'$,
$H_2=\{ z_2=0\}' = \{ w_2=0\}$ and $H_3=\sigma_3^{-1}
(C_2) = \{ w_1=0\}$; since $H_1\not\ni d$, $E(d)=\{ H_2,H_3\}$.
We have $\nu_1(d)=2=\nu_1(a)$.
Therefore, $E^1(d)=\emptyset$, $s_1(d)=0$ and $\cE_1(d)=E(d)$.
$\cF_1(d)=\cG_1(d)=\{ (g_3,2)\}$ provides a codimension
$0$ presentation of $\inv_1$ at $d$, and we get a codimension
$1$ presentation by taking
\[
N_1(d) = \{ w_3=0\},\qquad
\cH_1(d) = \{ (w_1 w_2^4 ,2)\} .
\]
Therefore, $\mu_2(c)=\frac{5}{2}$ and 
$D_2(d)=w_1^{\frac{1}{2}} w_2^2$, so that $\nu_2(d)=0$ and
\[
\inv_X (d) = (2,0,0) = \inv_X(c) !
\]
However,
\[
\mu_2(d) = \frac{5}{2} < \frac{7}{2} = \mu_2(c) ;
\]
i.e., $1\le\mu_X(d) < \mu_X(c)$, where $\mu_X=\mu_2$
(cf. (2.8) ff.).
We get a codimension $1$ presentation of $\inv_X=\inv_{1\frac{1}{2}}$
at $d$ by taking
\[
N_1(d) = \{ w_3=0\},\qquad
\cG_2(d) = \big\{ \big( D_2(d),1\big)\big\} .
\]
Therefore,
\[
S_{\inv_X}(d) = S_{\inv_1}(d) = \{ w_2 = w_3 = 0\} ,
\]
so we let $\sigma_4$ be the blowing-up with centre
$C_3=\{ w_2 = w_3 = 0\}$.
Then $\sigma_4^{-1}(U_{121}) = U_{1212} \cup U_{1213}$,
where $U_{121i}=\sigma_4^{-1} (U_{121})\bs \{ w_i = 0\}'$,
$i=2,3$;
$U_{1212}$ has coordinates $(v_1,v_2,v_3)$ in which $\sigma_4$
is given by
\[
w_1=v_1,\qquad
w_2=v_2,\qquad
w_3=v_2 v_3 .
\]

\medskip

{\em Year four.}\quad
Let $X_4$ be the strict transform of $X_3$.
Then $X_4\cap U_{1212}=V(g_4)$, where $g_4(v)=v_3^2-v_1 v_2^2$.
Let $e=0$ in $U_{1212}$.
Then $E(e)=\{ H_3,H_4\}$, where $H_3=\{ w_1=0\}'=\{v_1=0\}$
and $H_4=\sigma_4^{-1} (C_3)=\{v_2=0\}$.
Again $\nu_1(e)=2=\nu_1(a)$, so that $E^1(e)=\emptyset$,
$s_1(e)=0$ and $\cE_1(e)=E(e)$.
Calculating as above, we obtain $\mu_2(e)=\frac{3}{2}$ and
$D_2(e)=v_1^{\frac{1}{2}} v_2$, so that $\nu_2(e)=0$ and
$\inv_X(e)=(2,0;\, 0)$ again.
But now $\mu_X(e)=\mu_2(e)=3/2$.
Our invariant $\inv_X$ is presented at $e$ by
\[
N_1(e) = \{ v_3=0\},\qquad
\cG_2(e) = \{ (v_1^{\frac{1}{2}} v_2,1)\} .
\]
Therefore, $S_{\inv_X}(e)=\{ v_2=v_3=0\}$.
Taking as $\sigma_5$ the blowing-up with centre
$C_4=S_{\inv_X}(e)$, the strict transform $X_5$ becomes
smooth (over $U_{1212}$).
($\mu_2(e)-1<1$, so $\nu_1(\cdot)$ must decrease
over $C_4$.)

Further blowings-up are still needed to obtain the
stronger assertion of embedded resolution of singularities.

\rem {\em 3.7.}\quad
The hypersurface $V(g_4)$ in year four above is
called ``Whitney's umbrella''.
Consider the same hypersurface $X=\{x_3^2-x_1x_2^2=0\}$ but
without a history of blowings-up; i.e., $E(\cdot)=\emptyset$.
Let $a=0$.
In this case, $\inv_{1\frac{1}{2}}(a)=(2,0;\, \frac{3}{2})$, and
we get a codimension $1$ presentation of $\inv_{1\frac{1}{2}}$
at $a$ using
\[
N_1(a) =\{ x_3=0\},\qquad
\cG_2(a) = \{ (x_1 x_2^2 ,3)\}
\]
or, equivalently, $\cG_2(a)=\{ (x_1,1),(x_2,1)\}$, as in
year zero of Example 3.1.
Therefore,
\[
\inv_X(a) = (2,0;\, \frac{3}{2},0;\, 1,0;\, \infty) .
\]
As a centre of blowing up we would choose $C=S_{\inv_X}(a)=
\{ a\}$ --- not the $x_1$-axis as in year four above,
although the singularity is the same!

\section{Key properties of the invariant}

Our main goal in this section is to explain why $\inv_X(a)$
is indeed an invariant.
Once we establish invariance, the Embedded Desingularization
Theorem A follows directly from local properties of $\inv_X$.
The crucial properties have already been explained
in Section 3 above; we summarize them in the following theorem.

\theorem {\bf B.}\quad {\rm ([BM5, Th. 1.14].)}
{\em Consider any sequence of $\inv_X$-admissible (local)
blowings-up (1.8).
Then the following properties hold:

{\rm (1) Semicontinuity.} (i) For each $j$, every point
of $M_j$ admits a neighbourhood $U$ such that $\inv_X$
takes only finitely many values in $U$ and, for all $a\in U$,
$\{ x\in U:\ \inv_X(x)\le\inv_X(a)\}$ is Zariski-open
in $U$.
(ii) $\inv_X$ is {\rm infinitesimally upper-semicontinuous}
in the sense that $\inv_X(a)\le\inv_X\big(\sigma_j(a)\big)$
for all $a\in M_j$, $j\ge 1$.

{\rm (2) Stabilization.} Given $a_j\in M_j$ such that
$a_j=\sigma_{j+1}(a_{j+1})$, $j=0,1,2,\ldots$, there exists
$j_0$ such that $\inv_X(a_j)=\inv_X(a_{j+1})$ when $j\ge j_0$.
(In fact, any nonincreasing sequence in the value set
of $\inv_X$ stabilizes.)

{\rm (3) Normal crossings.}  Let $a\in M_j$.
Then $S_{\inv_X}(a)$ and $E(a)$ simultaneously have only
normal crossings.
Suppose $\inv_X(a)=\big(\ldots; \nu_{t+1} (a)\big)$.
If $\nu_{t+1}(a)=\infty$, then $S_{\inv_X}(a)$ is smooth.
If $\nu_{t+1}(a)=0$ and $Z$ denotes any irreducible
component of $S_{\inv_X}(a)$, then
\[
Z=S_{\inv_X} (a) \cap \bigcap \{ H\in E(a):\ Z\subset H\} .
\]

{\rm (4) Decrease.}  Let $a\in M_j$ and suppose
$\inv_X(a)=\big(\ldots; \nu_{t+1}(a)\big)$.
If $\nu_{t+1}(a)=\infty$ and $\sigma$ is the local
blowing-up of $M_j$ with centre $S_{\inv_X}(a)$, then
$\inv_X(a')<\inv_X(a)$ for all $a'\in\sigma^{-1}(a)$.
If $\nu_{t+1}(a)=0$, then there is an additional
invariant $\mu_X(a)=\mu_{t+1}(a)\ge 1$ such that,
if $Z$ is an irreducible component of $S_{\inv_X}(a)$
and $\sigma$ is the local blowing-up with centre $Z$,
then $\big( \inv_X(a'),\mu_X(a')\big)< \big( \inv_X(a),
\mu_X(a)\big)$ for all $a'\in\sigma^{-1}(a)$.
($e_t! \mu_X(a)\in\BN$, where $e_t$ is defined
as in Section 1 or in the proof following.)}

\prf
The semicontinuity property (1)(i) has been explained
in Remark 3.5.
Infinitesimal upper-semicontinuity (1)(ii) is immediate
from the definition of the $s_r(a)$ and from
infinitesimal upper-semicontinuity of the order
of a function on blowing up locally with 
smooth centre in its equimultiple locus.
(The latter property is an elementary Taylor series
computation, and is also clear from the calculation
in Section 2 above.)

The stabilization property (2) for $\inv_{\frac{1}{2}}$ is
obvious in the hypersurface case because then $\inv_{\frac{1}{2}} (a)=
\nu_1(a)\in\BN$.
(In the general case, we need to begin with stabiization
of the Hilbert-Samuel function;
see [BM2, Th. 5.2.1] for a very simple proof of this result
due originally to Bennett [Be].)
The stabilization assertion for $\inv_X$ follows from
that for $\inv_{\frac{1}{2}}$ and from infinitesimal
semicontinuity because, although $\nu_{r+1}(a)$, for
each $r>0$, is perhaps only rational, our construction
in Section 3 shows that $e_r! \nu_{r+1} (a)\in\BN$,
where $e_1=\nu_1(a)$ and $e_{r+1}=\max \{ e_r !,
e_r! \nu_{r+1} (a)\}$, $r>0$.
(In the general case, the Hilbert-Samuel function
$H_{X_j,a}(\ell)$ coincides with a polynomial if
$\ell\ge k$, for $k$ large enough, and we can take
as $e_1$ the least such $k$.)

The normal crossings condition (3) has also been explained
in Section 3; see Remark 3.6, in particular, for the case
that $\nu_{t+1}(a)=0$.
The calculation in Section 2 then gives the property
of decrease (4), as is evident also in the example
of Section 3.
\qed

\medskip

When our spaces satisfy a compactness assumption
(so that $\inv_X$ takes maximum values), it follows
from Theorem B that we can obtain the Embedded
Desingularization Theorem A by simply applying the
algorithm of 1.11 above, stopping when $\inv_X$
becomes (locally) constant.
To be more precise, let $\inv_X^\rme$ denote the
extended invariant for canonical desingularization
introduced in Remark 3.6.
Consider a sequence of blowings-up (1.8) with
$\inv_X$-admissible centres.
Note that if $X_j$ is not smooth and
$a\in \Sing\, X_j$, then $S_{\inv_X}(a)\subset
\Sing\, X_j$ because $\nu_1$ (or, in general,
$H_{X_j,a}$) already distinguishes between smooth
and singular points.
Since $\Sing\, X_j$ is Zariski-closed, it follows
that if $C_j$ denotes the locus of maximum values
of $\inv_X^\rme$ on $\Sing\, X_j$, then $C_j$ is smooth.
By Theorem B, there is a finite sequence of blowings-up
with such centres, after which $X_j$ is smooth.

On the other hand, if $X_j$ is smooth and $a\in S_j$,
where $S_j=\{ x\in X_j:\ s_1(x)>0\}$, then $S_{\inv_X}(a)
\subset S_j$.
Since $S_j$ is Zariski-closed, it follows that if $C_j$
denotes the locus of maximum values of $\inv_X^\rme$ on
$S_j$, then $C_j$ is smooth.
Therefore, after finitely many further blowings-up
$\sigma_{j+1},\ldots,\sigma_k$ with such centres,
$S_k=\emptyset$.
It is clear from the definition of $s_1$ that,
if $X_k$ is smooth and $S_k=\emptyset$, then each
$H\in E_k$ which intersects $X_k$ is the strict
transform in $M_k$ of $\sigma_{i+1}^{-1} (C_i)$,
for some $i$ such that $X_i$ is smooth along
$C_i$; therefore, $X_k$ and $E_k$ simultaneously have
only normal crossings, and we have Theorem A.

\medskip

We will prove invariance of $\inv_X$ using the idea
of a ``presentation'' introduced in Section 3 above.
It will be convenient to consider ``presentation'' in
an abstract sense, rather than associated to a particular
invariant: Let $M$ denote a manifold and let $a\in M$.

\medskip

\noindent {\em Definitions 4.1.}\quad 
An abstract {\em (infinitesimal) presentation} of {\em codimension}
$p$ at $a$ means simply a triple ($N=N_p(a)$, $\cH(a)$, $\cE(a)$)
as in Section 3; namely: $N$ is a germ of a submanifold
of codimension $p$ at $a$, $\cH(a)$ is a finite collection
of pairs $(h,\mu_h)$, where $h\in\cO_{N,a}$, $\mu_h\in\BQ$ 
and
$\mu_a(h)\ge \mu_h$, and $\cE(a)$ is a finite set of smooth
hypersurfaces containing $a$, such that $N$ and $\cE(a)$
simultaneously have normal crossings and $N\not\subset H$,
for all $H\in\cE(a)$.

A local blowing-up $\sigma$ with centre $C\ni a$ will be
called {\em admissible} (for an infinitesimal presentation
as above) if $C\subset S_{\cH(a)}=\{ x\in N:\ \mu_x (h)\ge
\mu_h$, for all $(h,\mu_h)\in\cH(a)\}$.

\defn {\em 4.2.}\quad 
We will say that two infinitesimal presentations
($N=N_p(a)$, $\cH(a)$, $\cE(a)$) and ($P=P_q(a)$,
$\cF(a)$, $\cE(a)$) with given $\cE(a)$, but not necessarily
of the same codimension, are {\em equivalent} if
(in analogy with Definition 3.2):

(1) $S_{\cH(a)}=S_{\cF(a)}$, as germs at $a$ in $M$.

(2) If $\sigma$ is an admissible local blowing-up and
$a'\in\sigma^{-1}(a)$, then $a'\in N'$ and $\mu_{a'}
(y_\exc^{-\mu_h} h\circ\sigma)\ge\mu_h$ for all
$(h,\mu_h)\in \cH(a)$ if and only if $a'\in P'$ and
$\mu_{a'}(y_\exc^{-\mu_f} f\circ\sigma)\ge\mu_f$
for all $(f,\mu_f)\in\cF(a)$.

(3) Conditions (1) and (2) continue to hold for the
transforms ($N_p(a')$, $\cH(a')$, $\cE(a')$) and
($P_q(a')$, $\cF(a')$, $\cE(a')$) of our data by
sequences of morphisms of types (i), (ii) and (iii) as
in Definition 3.2.

\medskip

We will, in fact, impose a further condition on the
way that exceptional blowings-up (iii) are allowed
to occur in a sequence of transformations in condition
(3) above; see Definition 4.5 below.

Our proof of invariance of $\inv_X$ follows the
constructive definition outlined in Section 3.
Let $X$ denote a hypersurface in $M$, and consider
any sequence of blowings-up (or local blowings-up)
(1.8), where we assume (at first) that the centres
of blowing up are $\frac{1}{2}$-admissible.
Let $a\in M_j$, for some $j=0,1,2,\ldots$.
Suppose that $g\in \cO_{M_j,a}$ generates the local ideal
$\cI_{X_j,a}$ of $X_j$ at $a$, and let $\mu_g=\mu_a(g)$.
Then, as in Section 3, $\cG_1(a)=\{ (g,\mu_g)\}$ determines
a codimension zero presentation ($N_0(a)$, $\cG_1(a)$,
$\cE_0(a)$) of $\inv_{\frac{1}{2}}=\nu_1$ at $a$, where
$N_0(a)$ is the germ of $M_j$ at $a$, and $\cE_0(a)=\emptyset$.
In particular, the equivalence class of ($N_0(a)$,
$\cG_1(a)$, $\cE_0(a)$) in the sense of Definition 4.2
depends only on the local isomorphism class of $(M_j,X_j)$
at $a$.

We introduce $E^1(a)$ as in 1.12 above, and let
$s_1(a)=\# E^1(a)$, $\cE_1(a)=E(a)\backslash E^1(a)$.
Let
\[
\cF_1(a) = \cG_1(a) \cup \big( E^1 (a),1\big) ,
\]
where $\big( E^1(a),1\big)$ denotes $\{ (x_H,1):\ H\in
E^1(a)\}$ and $x_H$ means a local generator of the ideal
of $H$.
Then ($N_0(a)$, $\cF_1(a)$, $\cE_1(a)$) is a codimension
zero presentation of $\inv_1=(\nu_1,s_1)$ at $a$.
Clearly, the equivalence class of ($N_0(a)$, $\cF_1(a)$,
$\cE_1(a)$) depends only on the local isomorphism class of
($M_j$, $X_j$, $E_j$, $E^1(a)$).
Moreover, ($N_0(a)$, $\cF_1(a)$, $\cE_1(a)$) has an equivalent
codimension one presentation ($N_1(a)$, $\cH_1(a)$, $\cE_1(a)$)
as described in Section 3.
For example, let $a_k=\pi_{kj} (a)$, $k=0,\ldots,j$,
as in 1.12, and let $i$ denote the ``earliest year''
$k$ such that $\inv_{\frac{1}{2}}(a)=\inv_{\frac{1}{2}}(a_k)$.
Then $\cE_1(a_i)=\emptyset$.
As in Section 3, we can take $N_1(a_i) =$ any hypersurface
of maximal contact for $X_i$ at $a_i$.
If $(x_1,\ldots,x_n)$ are local coordinates for $M_i$
with respect to which $N_1(a_i)=\{ x_n=0\}$, then we
can take
\[
\cH_1(a_i) = \left\{ \left( \frac{\partial^q f}{\partial x_n^q}
\bigg|_{N_1(a_i)},\, \mu_f-q\right):\ 0\le q<\mu_f,\ (f,\mu_f)
\in \cF_1(a_i)\right\} .
\]
A codimension one presentation ($N_1(a)$, $\cH_1(a)$,
$\cE_1(a)$) of $\inv_1$ at $a$ can be obtained by transforming
($N_1(a_i)$, $\cH_1(a_i)$, $\cE_1(a_i)$) to $a$.
The condition that $N_1(a)$ and $\cE_1(a)$ simultaneously
have normal crossings and $N_1(a)\not\subset H$ for all
$H\in \cE_1(a)$ is a consequence of the effect of
blowing with smooth centre of codimension at least $1$
in $N(a_k)$, $i\le k<j$ (as in the calculation in
Section 2).

Say that $\cH_1(a)=\{ (h,\mu_h)\}$; each $h\in \cO_{N_1(a),a}$
and $\mu_h\le \mu_a(h)$.
Recall that we define
\begin{eqnarray*}
\mu_2(a) & = & \min_{\cH_1(a)} \frac{\mu_a(h)}{\mu_h}\\
\mu_{2H} (a)  & = & \min_{\cH_1(a)} \frac{\mu_{H,a}(h)}{\mu_h} ,
\qquad H \in \cE_1 (a) ,\\
\hbox{and}\qquad \nu_2(a) & = &\mu_2(a) - \sum_{H\in\cE_1(a)}
\mu_{2H}(a) .
\end{eqnarray*}
(Definitions 3.2, 3.4).
Propositions 4.4 and 4.6 below show that each of $\mu_2(a)$
and $\mu_{2H}(a)$, $H\in\cE_1(a)$, depends only on the
equivalence class of ($N_1(a)$, $\cH_1(a)$, $\cE_1(a)$),
and thus only on the local isomorphism class of ($M_j$,
$X_j$, $E_j$, $E^1(a)$).

If $\nu_2(a)=0$ or $\infty$, then we set $\inv_X(a)=\inv_
{1\frac12}(a)$.
If $0<\nu_2(a)<\infty$, then we construct a codimension
one presentation ($N_1(a)$, $\cG_2(a)$, $\cE_1(a)$) of
$\inv_{1\frac12}$ at $a$, as in Section 3.
From the construction, it is not hard to see that the
equivalence class of ($N_1(a)$, $\cG_2(a)$, $\cE_1(a)$) depends
only on that of ($N_1(a)$, $\cH_1(a)$, $\cE_1(a)$).
(See [BM5, 4.23 and 4.24] as well as Proposition 4.6 ff. below.)

This completes a cycle in the inductive definition of $\inv_X$.
Assume now that the centres of the blowings-up in (1.8)
are $1\frac12$-admissible.
We introduce $E^2(a)$ as in 1.12, and let $s_2(a)=\#E^2(a)$,
$\cE_2(a)=\cE_1(a)\backslash E^2(a)$.
If $\cF_2(a)=\cG_2(a)\cup \big( E^2(a),1\big)$, where
$\big( E^2(a),1\big)$ denotes $\{ (x_H|_{N_1(a)},1):\ 
H\in E_2(a)\}$, then ($N_1(a), \cF_2(a), \cE_2(a)$) is a
codimension one presentation of $\inv_2 = (\inv_{1\frac12},s_2)$
at $a$, whose equivalence class depends only on the local
isomorphism class of ($M_j$, $X_j$, $E_j$, $E^1(a)$, $E^2(a)$).
It is clear from the construction of $\cG_2(a)$ that
$\mu_{\cG_2(a)}=1$, where
\[
\mu_{\cG_2(a)} = \min_{(g,\mu_g)\in\cG_2(a)}
\frac{\mu_a(g)}{\mu_g} .
\]
Therefore ($N_1(a)$, $\cF_2(a)$, $\cE_2(a)$) admits an
equivalent codimension two presentation ($N_2(a)$, $\cH_2(a)$,
$\cE_2(a)$), and we define $\nu_3(a)=\mu_3(a)-\sum_{H\in\cE_2(a)}
\mu_{3H}(a)$, as above.
By Propositions 4.4 and 4.6,
 $\mu_3(a)$ and each $\mu_{3H}(a)$ depend only on the
equivalence class of ($N_2(a)$, $\cH_2(a)$, $\cE_2(a)$), $\ldots$.
We continue until $\nu_{t+1}(a)=0$ or $\infty$ for some $t$,
and then take $\inv_X(a)=\inv_{t+\frac12}(a)$.

Invariance of $\inv_X$ thus follows from Propositions 4.4 and
4.6 below, which are formulated purely in terms of an abstract
infinitesimal presentation.

Let $M$ be a manifold, and let ($N(a)$, $\cH(a)$, $\cE(a)$)
be an infinitesimal presentation of codimension $r\ge 0$
at a point $a\in M$.
We write $\cH(a)=\{ (h,\mu_h)\}$, where $\mu_a(h)\ge\mu_h$
for all $(h,\mu_h)$.

\medskip

\noindent {\em Definitions 4.3.}\quad 
We define $\mu(a)=\mu_{\cH(a)}$ as
\[
\mu_{\cH(a)} = \min_{\cH(a)} \frac{\mu_a(h)}{\mu_h} .
\]
Thus $1\le \mu(a)\le\infty$.
If $\mu(a)<\infty$, then we define $\mu_H(a)=\mu_{\cH(a),H}$,
for each $H\in\cE(a)$, as
\[
\mu_{\cH(a),H} = \min_{\cH(a)} \frac{\mu_{H,a}(h)}{\mu_h} .
\]

\medskip

We will show that each of $\mu(a)$ and the $\mu_H(a)$ depends
only on the equivalence class of ($N(a)$, $\cH(a)$, $\cE(a)$)
(where we consider only {\em presentations of the same
codimension} $r$).
The main point is that $\mu(a)$ and the $\mu_H(a)$ can be
detected by ``test blowings-up'' (test transformations
of the form (i), (ii), (iii) as allowed by the definition
4.2 of equivalence).

For $\mu(a)$, we show in fact that if ($N^i(a)$,
$\cH^i(a)$, $\cE(a)$), $i=1,2$, are two infinitesimal
presentations of the same codimension $r$, then
$\mu_{\cH^1(a)}=\mu_{\cH^2(a)}$ if the presentations
are equivalent with respect to transformations of
types (i) and (ii) alone (i.e., where we allow only
transformations of types (i) and (ii) in
Definition 4.2).
This ia a stronger condition than invariance under
equivalence in the sense of Definition 4.2 (using all
three types of transformations) because the equivalence
class with respect to transformations of types (i) and
(ii) alone is, of course, larger than the equivalence
class with respect to transformations of all three
types (i), (ii) and (iii).

\prop {\bf 4.4.} 
[BM5, Prop. 4.8].
{\em $\mu(a)$ depends only on the equivalence class of
($N(a)$, $\cH(a)$, $\cE(a)$) (among presentations
of the same codimension $r$) with respect to
transformations of types (i) and (ii).}

\prf
Clearly, $\mu(a)=\infty$ if and only if $S_{\cH(a)}=N(a)$;
i.e., if and only if $S_{\cH(a)}$ is (a germ of) a
submanifold of codimension $r$ in $M$.

Suppose that $\mu(a)<\infty$.
We can assume that $\cH(a)=\{ (h,\mu_h)\}$ where all $\mu_h=e$,
for some $e\in \BN$.
Let $\sigma_0$: $P_0=M\times\BK\to M$ be the projection
from the product with a line (i.e., a morphism of
type (ii)) and let ($N(c_0)$, $\cH(c_0)$, $\cE(c_0)$) denote
the transform of ($N(a)$, $\cH(a)$, $\cE(a)$) at $c_0=
(a,0)\in P_0$; i.e., $N(c_0)=N(a)\times\BK$, $\cE(c_0)=
\{ H\times\BK$, for all $H\in\cE(a)$, and $M\times\{0\}\}$ and
$\cH(c_0)=\{ (h\circ\sigma_0,\mu_h)$: $(h,\mu_h)\in\cH(a)\}$.
We follow $\sigma_0$ by a sequence of admissible
blowings-up (morphisms of type (i)),
\[
\lra\ P_{\beta+1}\ 
\stackrel{\sigma_{\beta+1}}{\lra}\ 
P_\beta\ \lra\ \cdots\ \lra\ P_1\ 
\stackrel{\sigma_1}{\lra}\ P_0 ,
\]
where each $\sigma_{\beta+1}$ is a blowing-up with centre a
point $c_\beta\in P_\beta$ determined as follows:
Let $\gamma_0$ denote the arc in $P_0$ given by
$\gamma_0(t)=(a,t)$.
For $\beta\ge 1$, define $\gamma_{\beta+1}$ inductively
as the lifting of $\gamma_\beta$ to $P_{\beta+1}$, and set
$c_{\beta+1}=\gamma_{\beta+1}(0)$.

We can choose local coordinates $(x_1,\ldots,x_n)$ for $M$
at $a$, in which $a=0$ and $N(a)=\{ x_{n-r+1}=\cdots=
x_n=0\}$.
Write $(x,t)=(x_1,\ldots,x_{n-r},t)$ for the corresponding
coordinate system of $N(c_0)$.
In $P_1$, the strict transform $N(c_1)$ of $N(c_0)$ has a
local coordinate system $(x,t)=(x_1,\ldots,x_{n-r},t)$
at $c_1$ with respect to which $\sigma_1(x,t)=(tx,t)$,
and $\gamma_1(t)=(0,t)$ in this coordinate chart; moreover,
$\cH(c_1)=\{ (t^{-e}h(tx),e)$, for all $(h,\mu_h)=
(h,e)\in\cH(a)\}$.
After $\beta$ blowings-up as above, $N(c_\beta)$ has a local
coordinate system $(x,t)=(x_1,\ldots,x_{n-r},t)$ with respect
to which $\sigma_1\circ\cdots\circ\sigma_\beta$ is given
by $(x,t)\mapsto (t^\beta x,t)$, $\gamma_\beta (t)=(0,t)$
and $\cH(c_\beta)=\{ (h',\mu_{h'}=e)\}$, where
\[
h' = t^{-\beta e} h(t^\beta x) ,
\]
for all $(h,\mu_h)=(h,e)\in\cH(a)$.
By the definition of $\mu(a)$, each
\[
h(t^\beta x) = t^{\beta \mu(a)e} \th'(x,t) ,
\]
where the $\th'(x,t)$ do not admit $t$ as a common
divisor; for each $(h,\mu_h)\in\cH(a)$, we have
$$
h' = t^{\beta (\mu(a)-1) e} \th' .
$$

We now introduce a subset $S$ of $\BN\times\BN$ depending
only on the equivalence class of ($N(a)$, $\cH(a)$,
$\cE(a)$) (with respect to transformations of types (i)
and (ii)) as follows:
First, we say that $(\beta,0)\in S$, $\beta\ge 1$, if
after $\beta$ blowings-up as above, there exists
(a germ of) a submanifold $W_0$ of codimension $r$ in
the exceptional hypersurface $H_\beta = \sigma_\beta^{-1}
(c_{\beta-1})$ such that $W_0\subset S_{\cH(c_\beta)}$.
If so, then necessarily $W_0=H_\beta \cap N(c_\beta)=\{t=0\}$,
and the condition that $W_0\subset S_{\cH(c_\beta)}$ means
precisely that $\mu_{W_0,c_\beta}(h')\ge e$, for all $h'$;
i.e., that $\beta\big(\mu(a)-1\big)e\ge e$, or
$\beta\big(\mu(a)-1\big)\ge 1$.
(In particular, since $\mu(a)\ge 1$, $(\beta,0)\not\in S$ for
all $\beta\ge 1$ if and only if $\mu(a)=1$.)

Suppose that $(\beta,0)\in S$, for some $\beta\ge 1$, as
above.
Then we can blow up $P_\beta$ locally with centre $W_0$.
Set $Q_0=P_\beta$, $d_0=c_\beta$ and $\delta_0=\gamma_\beta$.
Let $\tau_1$: $Q_1\to Q_0$ denote the local blowing-up
with centre $W_0$, and let $d_1=\delta_1(0)$, where
$\delta_1$ denotes the lifting of $\delta_0$ to $Q_1$.
(Then $\tau_1|N(d_1)$: $N(d_1)\to N(d_0)$ is the
identity.)
We say that $(\beta,1)\in S$ if there exists a submanifold
$W_1$ of codimension $r$ in the hypersurface $H_1=\tau_1^{-1} (W_0)$
such that $W_1\subset S_{\cH(d_1)}$.
If so, then again necessarily $W_0=H_1\cap N(d_1)=\{ t=0\}$.
Since $\cH(d_1)=\{ (h',e)\}$, where each $h'=t^{\beta
(\mu(a)-1)e-e}\th'$ and the $\th'$ do not admit $t$ 
as a common factor, it
follows that $(\beta,1)\in S$ if and only if $\beta\big(\mu(a)-1\big)
e-e\ge e$.

We continue inductively:
If $\alpha\ge 1$ and $(\beta,\alpha-1)\in S$, let $\tau_\alpha$:
$Q_\alpha\to Q_{\alpha-1}$ denote the local blowing-up with centre
$W_{\alpha-1}$, and let $d_\alpha=\delta_\alpha (0)$, where
$\delta_\alpha$ is the lifting of $\delta_{\alpha-1}$ to $Q_\alpha$.
We say that $(\beta,\alpha)\in S$ if there exists (a germ of) a
submanifold $W_\alpha$ of codimension $r$ in the exceptional
hypersurface $H_\alpha=\tau_\alpha^{-1} (W_{\alpha-1})$ such that
$W_\alpha\subset S_{\cH(d_\alpha)}$.
Since $\cH(d_\alpha)=\{ (h',e)\}$, where each
$h'=t^{\beta(\mu(a)-1)e-\alpha e}\th'$ and the $\th'$ do not
admit $t$ as a common factor, it follows as before that
$(\beta,\alpha)\in S$ if and only if $\beta\big(\mu(a)-1\big)-\alpha
\ge 1$.

Now $S$, by its definition, depends only on the equivalence
class of ($N(a)$, $\cH(a)$, $\cE(a)$) (with respect to
transformations of types (i) and (ii)).
On the other hand, we have proved that $S=\emptyset$ if and only if
$\mu(a)=1$, and, if $S\ne\emptyset$, then
\[
S = \big\{ (\beta,\alpha)\in\BN\times\BN:\ \beta\big( \mu(a)-1\big)
-\alpha\ge 1\big\} .
\]
Our proposition follows since $\mu(a)$ is uniquely determined by
$S$; in the case that $S\ne \emptyset$,
$$
\mu(a) = 1 + \sup_{(\beta,\alpha)\in S} \frac{\alpha+1}{\beta} .
\eqno{\Box}
$$

\medskip

Suppose that $\mu(a)<\infty$.
Then we can also use test blowings-up to prove invariance
of $\mu_H(a)=\mu_{\cH(a),H}$, $H\in\cE(a)$:
Fix $H\in\cE(a)$.
As before we begin with the projection $\sigma_0$:
$P_0=M\times\BK\to M$ from the product with a line.
Let ($N(a_0)$, $\cH(a_0)$, $\cE(a_0)$) denote the
transform of ($N(a)$, $\cH(a)$, $\cE(a)$) at $a_0=
(a,0)\in P_0$ by the morphism $\sigma_0$ (of type (ii)),
and let $H_0^0 = M\times\{0\}$, $H_1^0 =\sigma_0^{-1}
(H) = H\times\BK$.
Thus $H_0^0,H_1^0\in \cE(a_0)$.
We follow $\sigma_0$ by a sequence of exceptional
blowings-up (morphisms of type (iii)),
\[
\lra\ P_{j+1}\ 
\stackrel{\sigma_{j+1}}{\lra}\ P_j\ 
\lra\ \cdots\ \lra\ P_1\ 
\stackrel{\sigma_1}{\lra}\ P_0 ,
\]
where each $\sigma_{j+1}$, $j\ge 0$, has centre
$C_j=H_0^j\cap H_1^j$ and $H_0^{j+1}=\sigma_{j+1}^{-1}
(C_j)$, $H_1^{j+1}=$ the strict transform of $H_1^j$ by
$\sigma_{j+1}$.
Let $a_{j+1}$ denote the unique intersection point of
$C_{j+1}$ and $\sigma_{j+1}^{-1} (a_j)$, $j\ge 0$.
($a_{j+1}=\gamma_{j+1}(0)$, where $\gamma_0$ denotes
the arc $\gamma_0(t)=(a,t)$ in $P_0$ and $\gamma_{j+1}$
denotes the lifting of $\gamma_j$ by $\sigma_{j+1}$,
$j\ge 0$.)

We can choose local coordinates $(x_1,\ldots,x_n)$ for $M$
at $a$, in which $a=0$, $N(a)=\{ x_{n-r+1} = \cdots = x_n=0\}$,
and each $K\in\cE(a)$ is given by $x_i=0$, for some
$i=1,\ldots,n-r$.
(Set $x_i=x_K$.)
Write $(x,t)=(x_1,\ldots,x_m,t)$, where $m=n-r$, for the
corresponding coordinate system of $N(a_0)=N(a)\times\BK$.

We can assume that $x_H=x_1$.
In $P_1$, the strict transform $N(a_1)$ of $N(a_0)$ has a chart
with coordinates $(x,t)=(x_1,\ldots,x_m,t)$ in which $\sigma_1$
is given by $\sigma_1(x,t)=(tx_1,x_2,\ldots,x_m,t)$ and in
which $a_1=(0,0)$, $\gamma_1(t)=(0,t)$ and $x_1=x_H$.
($x_H$ now means $x_{H_1^1}$.)
Proceeding inductively, for each $j$, $N(a_j)$ has a coordinate
system $(x,t)=(x_1,\ldots,x_m,t)$ in which $a_j=(0,0)$ and
$\sigma_1\circ\cdots\circ\sigma_j$: $N(a_j)\to N(a_0)$ is
given by
\[
(x,t)\mapsto (t^j x_1, x_2,\ldots,x_m,t) .
\]

We can assume that $\mu_h=e\in\BN$, for all
$(h,\mu_h)\in\cH(a)$.
Set
$$
D = \prod_{K\in\cE(a)} x_K ^{\mu_K(a)} .
$$
Thus $D^e$ is a monomial in the coordinates
$(x_1,\ldots,x_m)$ of $N(a)$ with exponents in $\BN$, and
$D^e$ is the greatest common divisor of the $h$ in
$\cH(a)$ which is a monomial in $x_K$, $K\in\cE(a)$ (by
Definitions 4.3).
In particular, for some $h=D^e g$ in $\cH(a)$, $g=g_H$ is
not divisible by $x_1=x_H$.
Therefore, there exists $i\ge 1$ such that
\[
\mu_{a_j} (g_H\circ\pi_j) = \mu_{a_i} (g_H\circ\pi_i) ,
\]
for all $j\ge i$, where $\pi_j=\sigma_0\circ\sigma_1\circ
\cdots\circ\sigma_j$.
(We can simply take $i$ to be the least order of a monomial
not involving $x_H$ in the Taylor expansion of $g_H$.)

On the other hand, for each $h=D^e g$ in $\cH(a)$,
$\mu_{a_j}(g\circ\pi_j)$ increases as $j\to\infty$ unless $g$
is not divisible by $x_H$.
Therefore, we can choose $h=D^e g_H$, as above, and $i$
large enough so that we also have $\mu(a_j)=\mu_{a_j}
(h\circ\pi_j)/e$, for all $j\ge i$.
Clearly, if $j\ge i$, then
\[
\mu_H (a) = \mu(a_{j+1}) - \mu(a_j) .
\]

Since $\mu(a)$ depends only on the equivalence class of
($N(a)$, $\cH(a)$, $\cE(a)$) among presentations of the
same codimension $r$, as defined by 4.2, the preceding
argument shows that each $\mu_H (a)$, $H\in\cE(a)$, is also
an invariant of this equivalence class.
But the argument shows more precisely that the $\mu_H(a)$
depend only on a larger equivalence class obtained by
allowing in Definition 4.2 only certain sequences
of morphisms of types (i), (ii) and (iii):

\defn {\em 4.5.}\quad 
We weaken the notion of equivalence in Definition 4.2
by allowing only the transforms induced by certain
sequences of morphisms of types (i), (ii) and (iii);
namely,
\[
\begin{array}{cccccccccccc}
\rightarrow & M_j &\stackrel{\sigma_j}{\rightarrow} &M_{j-1}
&\rightarrow &\cdots &\stackrel{\sigma_{i+1}}{\rightarrow}
& M_i &\rightarrow &\cdots &\rightarrow &M_0=M \\
&\cE(a_j)&&\cE(a_{j-1})&&&&\cE(a_i)&&&&\cE(a_0)=\cE(a)
\end{array}
\]
where, if $\sigma_{i+1},\ldots,\sigma_j$ are exceptional
blowings-up (iii), then $i\ge 1$ and $\sigma_i$ is of
either type (iii) or (ii).
In the latter case, $\sigma_i$: $M_i=M_{i-1}\times\BK\to
M_{i-1}$ is the projection, each $\sigma_{k+1}$,
$k=i,\ldots,j-1$, is the blowing-up with centre
$C_k=H_0^k \cap H_1^k$ where $H_0^k$, $H_1^k\in\cE(a_k)$,
$a_{k+1}=\sigma_{k+1}^{-1}(a_k)\cap H_1^{k+1}$, and we
require that the $H_0^k$, $H_1^k$ be determined by some
fixed $H\in \cE(a_{i-1})$ inductively in the following way:
$H_0^i=M_{i-1}\times\{0\}$, $H_1^i=\sigma_i^{-1} (H)$,
and, for $k=i+1,\ldots,j-1$, $H_0^k=\sigma_k^{-1}(C_{k-1})$,
$H_1^k=$ the strict transform of $H_1^{k-1}$ by $\sigma_k$.
\medskip

In other words, with this notion of equivalence, we have proved:

\prop {\bf 4.6.}\quad 
[BM5, Prop. 4.11].
{\em Each $\mu_H(a)$, $H\in\cE(a)$, and therefore also
$\nu(a)=\mu(a)-\Sigma \mu_H(a)$ depends only on the
equivalence class of ($N(a)$, $\cH(a)$, $\cE(a)$)
(among presentations of the same codimension).}
\medskip

Recall that in the $r$'th cycle of our recursive
definition of $\inv_X$, we use a codimension $r$ presentation
($N_r(a)$, $\cH_r(a)$, $\cE_r(a)$) of $\inv_r$ at $a$
to construct a codimension $r$ presentation
($N_r(a)$, $\cG_{r+1}(a)$, $\cE_r(a)$) of $\inv_{r+\frac12}$
at $a$.
The construction involved survives transformations as
allowed by Definition 4.5, but perhaps not an arbitrary
sequence of transformations of types (i), (ii) and (iii)
(cf. [BM5, 4.23 and 4.24]; in other words, we show only
that the equivalence class of ($N_r(a)$, $\cG_{r+1}(a)$,
$\cE_r(a)$) as given by Definition 4.5 depends only on that
of ($N_r(a)$, $\cH_r(a)$, $\cE_r(a)$).
It is for this reason that we need Proposition 4.6
as stated.

\vskip .25in

\noindent {\it Acknowledgement.}\quad We are happy to thank Paul Centore 
for the line drawings in this paper.

\end{document}